\documentclass[11pt,onecolumn]{article}

\usepackage{palatino}
\usepackage{mathpazo}

\usepackage[english]{babel}
\usepackage[latin1]{inputenc}
\usepackage{enumerate}
\usepackage{color}
\usepackage[T1]{fontenc}
\usepackage{epstopdf}
\usepackage{subfigure}
\usepackage{dsfont}
\usepackage[final]{graphicx}
\DeclareGraphicsExtensions{.eps}
\usepackage[T1]{fontenc}
\usepackage{amsmath}
\usepackage{mathtools}
\usepackage{amsthm}
\usepackage{amstext}
\usepackage{amssymb}
\usepackage{mathrsfs}
\usepackage{cite}
\usepackage{mathtools}
\usepackage{tikz}
\usetikzlibrary{arrows,positioning}

\def\R{\mathbb{R}}

\def\rP{\mathbb{P}}

\def\Gr{\mathop{\rm Gr}}

\def\B{{\mathcal B}}

\def\P{{\mathcal P}}

\def\M{{\mathcal M}}

\def\bmu{{\boldsymbol \mu}}
\def\bnu{{\boldsymbol \nu}}
\def\bpi{{\boldsymbol \pi}}
\def\bxi{{\boldsymbol \xi}}

\def\ba{{\bf a}}
\def\bx{{\bf x}}

\def\tpi{{\tilde \pi}}

\def\bs{{\bf s}}
\def\ts{{\tilde s}}
\def\hs{{\hat s}}
\def\bc{{\bf c}}

\def\sX{{\mathsf X}}

\def\sA{{\mathsf A}}
\def\sH{{\mathsf H}}
\def\sZ{{\mathsf Z}}
\def\sM{{\mathsf M}}
\def\sE{{\mathsf E}}

\def\sG{{\mathsf G}}

\def\sS{{\mathsf S}}

\theoremstyle{remark}

\newtheorem{definition}{Definition}
\newtheorem{theorem}{Theorem}
\newtheorem{corollary}{Corollary}

\newtheorem{proposition}{Proposition}
\newtheorem{lemma}{Lemma}
\theoremstyle{remark}
\newtheorem{remark}{Remark}
\newtheorem{assumption}{Assumption}
\allowdisplaybreaks

\newcommand{\appsec}{
\renewcommand{\thesubsection}{\Alph{subsection}}
}

\newcommand{\Halmos}{}

\begin{document}
\sloppy

\title{Discrete-time Risk-sensitive Mean-field Games}
\author{Naci Saldi, Tamer Ba\c{s}ar, and Maxim Raginsky
\thanks{N. Saldi is with the Department of Natural and Mathematical Sciences, Ozyegin University, Cekmekoy, Istanbul, Turkey, Email: naci.saldi@ozyegin.edu.tr T. Ba\c{s}ar and M. Raginsky are with the Coordinated Science Laboratory, University of Illinois at Urbana-Champaign, Urbana, IL 61801, Email: \{basar1,maxim\}@illinois.edu}}
\date{}
\maketitle

\abstract{
In this paper, we study a class of discrete-time mean-field games under the infinite-horizon risk-sensitive discounted-cost optimality criterion. Risk-sensitivity is introduced for each agent (player) via an exponential utility function. In this game model, each agent is coupled with the rest of the population through the empirical distribution of the states, which affects both the agent's individual cost and its state dynamics. Under mild assumptions, we establish the existence of a mean-field equilibrium in the infinite-population limit as the number of agents ($N$) goes to infinity, and then show that the policy obtained from the mean-field equilibrium constitutes an approximate Nash equilibrium when $N$ is sufficiently large.
}

\section{Introduction}
This paper deals with discrete-time mean-field games under the infinite-horizon risk-sensitive discounted-cost optimality criterion, where risk-sensitivity is introduced via an exponential utility function (see \cite{BaRi14,AvBrFe97,ChSo1987,HoMa72}). Mean-field game theory studies non-cooperative stochastic dynamic games with a large number of agents, who are coupled through their dynamics and cost functions via the mean-field term (i.e., the empirical distribution of their states). If the agents are identical (that is, they enter the game symmetrically), as the population goes to infinity, the mean-field term converges to a deterministic probability distribution by the law of large numbers, which results in the decoupling of agents from each other. Due to this decoupling, in the infinite-population regime, each agent is faced with a classical stochastic control problem subject to a constraint on the distribution of the state at each time. This constraint dictates that the state distribution should be consistent with the total population behavior; that is, at each time step, the law of the state must be the same as the limiting law of the mean-field term. Therefore, in the infinite-population case, the problem turns into one of establishing the existence of a policy and a state distribution flow such that this policy is an optimal solution of the stochastic control problem when the total population behavior is modeled by this state distribution flow and the resulting distribution of each agent's state is same as the state distribution flow when the generic agent applies this policy. This equilibrium --- also known as \emph{mean-field equilibrium} --- behavior between policy and state distribution flow is captured by the Nash certainty equivalence (NCE) principle (\cite{HuMaCa06}). Then, the solution of this single-agent decision problem provides an approximation to Nash equilibrium of games with large (but finite) population sizes.

The purpose of this paper is to establish, in discrete time, the existence of mean-field equilibrium for a general class of mean-field game models with risk-sensitive discounted-cost criteria and to show that the policy in the mean-field equilibrium constitutes an approximate Nash equilibrium for finite-agent games with sufficiently many agents.

Mean-field game theory has been introduced independently by \cite{HuMaCa06} and \cite{LaLi07} to establish approximate Nash equilibria for continuous-time non-cooperative differential games with a large number of identical agents. As roughly explained above, the key underlying idea is that, under the Nash certainty equivalence principle, the decentralized game problem can be reduced to a single-agent decision problem and the equilibrium solution of this decision problem constitutes an approximate Nash equilibrium for games with a sufficiently large number of agents. Characterization of this equilibrium solution leads to a Fokker-Planck (FP) equation that describes the total mass behavior and a Hamilton-Jacobi-Bellman (HJB) equation that describes the optimal solution corresponding to this total mass behavior. We refer the reader to \cite{HuCaMa07,TeZhBa14,Hua10,BeFrPh13,Ca11,CaDe13,GoSa14,MoBa16} for studies of continuous-time mean-field games with different models and cost functions, such as games with major-minor players, risk-sensitive games, games with Markov jump parameters, and LQG games.

In the literature, relatively few results are available on discrete-time mean-field games. Prior works have mostly studied the setup where the state space is discrete (finite or countable) and the agents are coupled only through their cost functions. \cite{GoMoSo10} consider a discrete-time mean-field game with a finite state space over a finite horizon. A discrete-time mean-field game with countable state-space is studied in \cite{AdJoWe15}, under an infinite-horizon discounted cost criterion. \cite{Bis15} considers the average-cost setting, where the state space is a $\sigma$-compact Polish space and the agents are coupled only through their cost functions. Discrete-time mean-field games with linear individual dynamics are studied in \cite{ElLiNi13,MoBa15,NoNa13,MoBa16-cdc}. We note that all existing works in the discrete-time setup deal with risk-neutral cost functions. Indeed, our work appears to be the first one that studies discrete-time mean-field games within the risk-sensitive framework.

Prior works on risk-sensitive mean-field games have mostly studied the continuous-time setup, and it is known that analyses of continuous-time and discrete-time setups are quite different, requiring different sets of tools. \cite{TeZhBa14} studies a continuous-time mean-field game with nonlinear individual dynamics and a risk-sensitive cost function. Tembine et al. characterize the mean-field equilibrium via coupled HJB and FP equations and explicit solutions to these equations are given when the individual state dynamics are linear. \cite{Tem15} considers a continuous-time mean-field game with nonlinear individual dynamics, where state dynamics have $L^p$-norm structure. Stochastic maximum principle is used to characterize the optimal solution of the problem. \cite{DjTe16} studies partially-observed version of the continuous-time risk-sensitive mean-field game. Djehiche and Tembine establish a stochastic maximum principle for the characterization of the mean-field equilibrium. \cite{MoBa17,MoBa15}  consider continuous-time risk-sensitive mean-field games with linear individual dynamics. First a generic risk-sensitive optimal control problem is solved which yields mean-field equilibrium and then it is shown that the policies in mean-field equilibrium lead to an approximate Nash equilibrium for games with a sufficiently large number of agents. It is also shown that this approximate Nash equilibrium is partially equivalent to the approximate Nash equilibrium of a certain robust mean-field game problem.

This paper deals with discrete-time risk-sensitive mean-field games with Polish state and action spaces and with bounded one-stage cost functions. The generic cost function of each player is exponential utility of an infinite-horizon discounted-cost function. In this game model, the agents are coupled through the mean-field term, which affects both the individual dynamics and one-stage cost functions. In the infinite-population case, a generic agent is faced with a non-homogenous risk-sensitive stochastic control problem under the NCE principle.
Due to NCE principle, the classical techniques used to study risk-sensitive stochastic control problems are not sufficient to establish the existence of a mean-field equilibrium.
To do this, we have to employ the fixed-point approach that is used to obtain equilibria in classical game problems, along with the dynamic programming recursion for risk-sensitive costs. In Section~\ref{sec2}, we formulate the finite-agent discrete-time risk-sensitive game problem of the mean-field type. In Section~\ref{sec3}, we introduce the infinite-population mean-field game and define mean-field equilibrium. In Section~\ref{sec3-1}, we prove the existence of a mean-field equilibrium. In Section~\ref{sec4} we establish that the mean-field equilibrium policies lead to an approximate Nash equilibrium for finite-agent games with sufficiently many agents. Section~\ref{conc} lists some directions for future research.

In \cite{SaBaRa17} we have studied the risk-neutral version of this problem under a similar set of assumptions on the system components. At the higher level, the techniques used in this paper to show the existence of a mean-field equilibrium and to prove that the policies in mean-field equilibrium provide an approximate Nash equilibrium for games with sufficiently many agents are almost the same as in \cite{SaBaRa17}; however there are some highly non-trivial technical differences between the analysis of risk-sensitive and risk-neutral cost functions. For instance, in the risk-neutral case, the cost function is in an additive form, and this results in a contractive dynamic programming operator. However, in the risk-sensitive case, the cost function is in a multiplicative form, and therefore, the dynamic programming operator is not contractive, which complicates the analysis.
Indeed, to establish the existence of a mean-field equilibrium, we need to first establish an optimality criterion for non-homogenous risk-sensitive stochastic control problems, which is absent in the literature. In the risk-neutral case, this optimality result has already been established. Moreover, as opposed to the risk-neutral case, we prove the existence of approximate Nash equilibrium by transforming the original problem into an equivalent one using a certain state-aggregation technique. This also brings in additional complexity to the analysis, which is unavoidable.

\smallskip

\noindent\textbf{Notation.} For a metric space $\sE$ with metric $d_{\sE}$, we let $C_b(\sE)$ denote the set of all bounded continuous real functions on $\sE$, and $\P(\sE)$ denote the set of all Borel probability measures on $\sE$. For any $\sE$-valued random element $x$, ${\cal L}(x)(\,\cdot\,) \in \P(\sE)$ denotes the distribution of $x$. A sequence $\{\mu_n\}$ of measures on $\sE$ is said to converge weakly to a measure $\mu$ if $\int_{\sE} g(e) \mu_n(de)\rightarrow\int_{\sE} g(e) \mu(de)$ for all $g \in C_b(\sE)$. For any $\nu \in \P(\sE)$ and measurable real function $g$ on $\sE$, we define $\nu(g) \coloneqq \int g d\nu$. For any subset $B$ of $\sE$, we let $\partial B$ and $B^c$ denote the boundary and complement of $B$, respectively. The notation $v\sim \nu$ means that the random element $v$ has distribution $\nu$.
Unless otherwise specified, the term ``measurable" will refer to Borel measurability.

\section{Finite Player Game Model}\label{sec2}

In this paper, we study a discrete-time $N$-agent infinite-horizon stochastic game of mean-field type with a Polish state space $\sX$ and a Polish action space $\sA$. The dynamics of each agent are governed by an initial state distribution $\mu_0 \in \P(\sX)$ and a state transition kernel $p : \sX \times \sA \times \P(\sX) \to \P(\sX)$ as follows. For every $t \in \{0,1,2,\ldots\}$ and every $i \in \{1,2,\ldots,N\}$, let $x^N_i(t) \in \sX$ and $a^N_i(t) \in \sA$ denote, respectively, the state and the action of Agent~$i$ at time $t$, and let
\begin{align}
e_t^{(N)}(\,\cdot\,) \coloneqq \frac{1}{N} \sum_{i=1}^N \delta_{x_i^N(t)}(\,\cdot\,) \in \P(\sX) \nonumber
\end{align}
denote the empirical distribution of the states at time $t$, where $\delta_x\in\P(\sX)$ is the Dirac measure at $x$. For $t=0$, $(x^N_1(0),\ldots,x^N_N(0))$ are independent and identically distributed according to $\mu_0$, and, for each $t \ge 0$, transition to the next-state configuration $(x^N_1(t+1),\ldots,x^N_N(t+1))$ is distributed according to the probability law
\begin{align}\label{eq:state_spec}
\prod^N_{i=1} p\big(dx^N_i(t+1)\big|x^N_i(t),a^N_i(t),e^{(N)}_t\big).
\end{align}
Here, the dynamics of each agent are weakly coupled through the mean-field term $e^{(N)}_t$. In the infinite population limit, $N \rightarrow \infty$, we see that, at each time step $t$, the empirical distribution $e_t^{(N)}$ converges to some deterministic probability measure $\mu_t$ by the law of large numbers, and so, this coupling vanishes.

The next step is to specify how the agents select their actions at each time step. To that end, we introduce the history spaces $\sH_0 = \sX \times \P(\sX)$ and $\sH_{t}=(\sX\times\sA\times\P(\sX))^{t}\times (\sX\times\P(\sX))$ for $t=1,2,\ldots$, all endowed with product Borel $\sigma$-algebras. We endow the set $\P(\sX)$ with the topology of weak convergence, which makes it a Polish space. A \emph{policy} for a generic agent is a sequence $\pi=\{\pi_{t}\}$ of stochastic kernels on $\sA$ given $\sH_{t}$; we say that such a policy is \emph{Markov} if each $\pi_t$ is a Markov kernel on $\sA$ given $\sX$. The set of all policies for Agent~$i$ is denoted by $\Pi_i$, and the subset consisting of all Markov policies by $\sM_i$. Furthermore, we let $\sM_i^c$ denote the set of all Markov policies for Agent~$i$ that are weakly continuous; that is, $\pi=\{\pi_t\}\in\sM_i^c$ if for all $t\geq0$, $\pi_t: \sX \rightarrow \P(\sA)$ is continuous when $\P(\sA)$ is endowed with the weak topology. Let
\begin{align}
{\bf \Pi}^{(N)} = \prod_{i=1}^N \Pi_i, \text{ } {\bf \sM}^{(N)} = \prod_{i=1}^N \sM_i, \text{ } \text{and} \text{ } {\bf \sM}^{(N,c)} = \prod_{i=1}^N \sM_i^c. \nonumber
\end{align}
We let ${\boldsymbol \pi}^{(N)} \coloneqq (\pi^1,\ldots,\pi^N)$, $\pi^i \in \Pi_i$, denote the $N$-tuple of policies for all the agents in the game, which will be referred to simply as a `policy.' Under such a policy, the actions of agents at each time $t \ge 0$ are generated according to the probability law
\begin{align}\label{eq:policy_spec}
\prod^N_{i=1} \pi^i_t\big(da^N_i(t)\big|h^N_i(t)\big),
\end{align}
where
\begin{align}
h^N_i(0) = (x^N_i(0),e^{(N)}_0) \text{ } \text{and} \text{ } h^N_i(t) = (h^N_i(t-1),x^N_i(t),a^N_i(t-1),e^{(N)}_t) \nonumber
\end{align}
for $t \ge 1$ are the local histories observed by Agent~$i$ at each time step. When ${\boldsymbol \pi}^{(N)} \in {\bf \sM}^{(N)}$, \eqref{eq:policy_spec} becomes
$$
\prod^N_{i=1}\pi^i_t(da^N_i(t)|x^N_i(t)).
$$
The stochastic transition kernels in \eqref{eq:state_spec} and \eqref{eq:policy_spec}, together with the initial state distribution $\mu_0$, uniquely determine the probability law of all the states and actions for all $i \in \{1,\ldots,N\}$ and all $t \geq 0$. We will denote expectations with respect to this probability law by $E^{{\boldsymbol \pi}^{(N)}}\big[\cdot\big]$.

We now define the optimality criterion for the problem. The \emph{one-stage cost} function for a generic agent is a measurable function $c : \sX \times \sA \times \P(\sX) \to [0,\infty)$. For Agent~$i$, we introduce "risk sensitivity" into the objective function of the agent. Accordingly, the agent's infinite-horizon risk-sensitive cost under the initial distribution $\mu_0$ and a policy ${\boldsymbol \pi}^{(N)} \in {\bf \Pi}^{(N)}$ is given by
\begin{align}
V_i^{(N)}({\boldsymbol \pi}^{(N)}) &= \frac{1}{\lambda} \log\biggl( E^{{\boldsymbol \pi}^{(N)}}\biggl[ e^{\lambda\sum_{t=0}^{\infty}\beta^{t}c(x_{i}^N(t),a_{i}^N(t),e^{(N)}_t)}\biggr]\biggr), \nonumber
\end{align}
where $\beta \in (0,1)$ is the discount factor and $\lambda > 0$ is the risk factor. By using Taylor series expansion of $V_i^{(N)}({\boldsymbol \pi}^{(N)})$ around $\lambda=0$, we note that to first order in $\lambda$,
\begin{align}
&V_i^{(N)}({\boldsymbol \pi}^{(N)}) \approx \biggl( E^{{\boldsymbol \pi}^{(N)}}\biggl[ \sum_{t=0}^{\infty}\beta^{t}c(x_{i}^N(t),a_{i}^N(t),e^{(N)}_t) \biggr] \nonumber \\
&\phantom{xxxxxxxxxxxxxxx}+ \frac{\lambda}{2} \mathrm{Var}^{{\boldsymbol \pi}^{(N)}}\biggl[ \sum_{t=0}^{\infty}\beta^{t}c(x_{i}^N(t),a_{i}^N(t),e^{(N)}_t) \biggr]\biggr), \nonumber
\end{align}
where $\mathrm{Var}[\,\cdot\,]$ denotes the variance of a random variable (see \cite[p. 107]{BaRi14}). Here, the second term in this expansion penalizes the variation of the cost around its mean value. Note that when $\lambda \rightarrow 0$, the cost function $V_i^{(N)}({\boldsymbol \pi}^{(N)})$ corresponds to the (risk-neutral) infinite-horizon discounted cost.

Since $\frac{1}{\lambda}\log(\cdot)$ is a strictly increasing function, without loss of generality, for the game we consider only the part with expectation:
\begin{align}
J_i^{(N)}({\boldsymbol \pi}^{(N)}) &= E^{{\boldsymbol \pi}^{(N)}}\biggl[ e^{\lambda\sum_{t=0}^{\infty}\beta^{t}c(x_{i}^N(t),a_{i}^N(t),e^{(N)}_t)}\biggr]. \nonumber
\end{align}
Next, we introduce the equilibrium solution for the game, which is Nash equilibrium:
\begin{definition}
A policy ${\boldsymbol \pi}^{(N*)}= (\pi^{1*},\ldots,\pi^{N*})$ constitutes a \emph{Nash equilibrium} for the $N$-player game, if
\begin{align}
J_i^{(N)}({\boldsymbol \pi}^{(N*)}) = \inf_{\pi^i \in \Pi_i} J_i^{(N)}({\boldsymbol \pi}^{(N*)}_{-i},\pi^i) \nonumber
\end{align}
for each $i=1,\ldots,N$, where ${\boldsymbol \pi}^{(N*)}_{-i} \coloneqq (\pi^{j*})_{j\neq i}$.
\end{definition}

As explained in detail in \cite{SaBaRa17}, for the risk-neutral game, there are two challenges to obtain Nash equilibria in mean-field games considered here: the (almost) decentralized nature of the information structure of the problem and the so-called \emph{curse of dimensionality} (the solution of the problem becomes intractable when the number of states/actions and agents is large). Therefore, it is of interest to find an approximate decentralized equilibrium with reduced complexity. To that end, we define the following.
\begin{definition}\label{def1}
A policy ${\boldsymbol \pi}^{(N*)} \in {\bf \sM}^{(N)}$ is a \emph{Markov-Nash equilibrium} if
\begin{align*}
J_i^{(N)}({\boldsymbol \pi}^{(N*)}) &= \inf_{\pi^i \in \sM_i} J_i^{(N)}({\boldsymbol \pi}^{(N*)}_{-i},\pi^i)
\end{align*}
for each $i=1,\ldots,N$, and an \emph{$\varepsilon$-Markov-Nash equilibrium} (for a given $\varepsilon > 0$) if
\begin{align*}
J_i^{(N)}({\boldsymbol \pi}^{(N*)}) &\leq \inf_{\pi^i \in \sM_i} J_i^{(N)}({\boldsymbol \pi}^{(N*)}_{-i},\pi^i) + \varepsilon
\end{align*}
for each $i=1,\ldots,N$.
\end{definition}

Note that, according to this definition, the agents can only use their current local state information $x^N_i(t)$ to design their policies. Indeed, in practical applications, due to memory and computational constraints, agents typically have access only to their current local state information. In addition, in the discrete-time mean field literature, it is common to establish the existence of approximate Nash equilibria with local (decentralized) information structures (see \cite{AdJoWe15} \cite{Bis15}).

In this paper, we prove, for the risk-sensitive game just formulated, the existence of the $\varepsilon$-Markov-Nash equilibria for games with sufficiently many agents using mean-field approach. We follow the steps in the risk-neutral version of the problem \cite{SaBaRa17}, with however some subtle differences in the risk-sensitive case. Accordingly, we first consider the infinite population limit and prove the existence of equilibrium. Then, for the finite-$N$ case, we show that, if each agent adopts the equilibrium policy in the infinite population limit, then the resulting policy will be an approximate Markov-Nash equilibrium for all sufficiently large $N$. It is important to note that, although the policy in the mean-field equilibrium is an approximate Markov-Nash equilibrium for the finite-agent game problem, it is indeed a true Nash equilibrium in the infinite population regime. This follows from the fact that the set of Markov policies is sufficiently rich for optimality in the limiting case, as each agent is faced with a single-agent decision problem.

\section{Mean-field games and mean-field equilibria}\label{sec3}

In this section, we introduce the infinite population limit of the original game introduced in the preceding section. This game is referred to as mean-field game and its equilibrium is called mean-field equilibrium. Mean-field games are not games in the classical sense: they are single-agent stochastic control problems whose state distribution at each time step should satisfy a certain consistency condition. In other words, we have a single agent and model the overall behavior of (a large population of) other agents (i.e., empirical distribution of the states) by an exogenous \textit{state-measure flow} $\bmu := (\mu_t)_{t \ge 0} \subset \P(\sX)$  with a given initial condition $\mu_0$. By the law of large numbers, this measure flow $\bmu$ should also be consistent with the state distributions of this single agent when the agent acts optimally. The precise mathematical description of the problem is given as follows.

The mean-field game model for a generic agent is specified by
\begin{align}
\bigl( \sX, \sA, p, c, \mu_0 \bigr), \nonumber
\end{align}
where, as before, $\sX$ and $\sA$ are the Polish state and action spaces, respectively. The stochastic kernel $p : \sX \times \sA \times \P(\sX) \to \P(\sX)$ denotes the transition probability. The measurable function $c: \sX \times \sA \times \P(\sX) \rightarrow [0,\infty)$ is the one-stage cost function and $\mu_0$ is the distribution of the initial state.

Define the history spaces $\sG_0 = \sX$ and $\sG_{t}=(\sX\times\sA)^{t}\times \sX$ for $t=1,2,\ldots$, all endowed with product Borel $\sigma$-algebras. A \emph{policy} is a sequence $\pi=\{\pi_{t}\}$ of stochastic kernels on $\sA$ given $\sG_{t}$. The set of all policies is denoted by $\Pi$. A \emph{Markov} policy is a sequence $\pi=\{\pi_{t}\}$ of stochastic kernels on $\sA$ given $\sX$. The set of Markov policies is denoted by $\sM$.

We define the set of all state-measure flows with a given initial condition $\mu_0$ as $\M \coloneqq \bigl\{\bmu \in \P(\sX)^{\infty}: \mu_0 \text{ is fixed}\bigr\}$. Given any measure flow $\bmu \in \M$, the probabilistic evolution of the states and actions of a generic agent are as follows
\begin{align}
x(0) &\sim \mu_0, \nonumber \\
x(t) &\sim p(\,\cdot\,|x(t-1),a(t-1),\mu_{t-1}), \text{ } t=1,2,\ldots \nonumber \\
a(t) &\sim \pi_t(\,\cdot\,|g(t)), \text{ } t=0,1,\ldots, \nonumber
\end{align}
where $g(t) \in \sG_t$ is the state-action history  up to time $t$. According to the Ionescu Tulcea Theorem (see, e.g., \cite{HeLa96}), an initial distribution $\mu_0$ on $\sX$, a policy $\pi$, and a state-measure flow $\bmu$ define a unique probability measure $P^{\pi}$ on $(\sX \times \sA)^{\infty}$. The expectation with respect to $P^{\pi}$ is denoted by $E^{\pi}$. A policy $\pi^{*} \in \Pi$ is said to be optimal for $\bmu$ if
\begin{align}
J_{\bmu}(\pi^{*}) = \inf_{\pi \in \Pi} J_{\bmu}(\pi), \nonumber
\end{align}
where the infinite-horizon risk-sensitive cost of policy $\pi$ with measure flow $\bmu$ is given by
\begin{align}
J_{\bmu}(\pi) &\coloneqq  E^{\pi}\biggl[ e^{\lambda\sum_{t=0}^{\infty} \beta^t c(x(t),a(t),\mu_t)} \biggr]. \nonumber
\end{align}
Using this definition, we first define the set-valued mapping
\begin{align}
\Psi : \M \rightarrow 2^{\Pi} \nonumber
\end{align}
as $\Psi({\boldsymbol \mu}) = \{\pi \in \Pi: \pi \text{ is optimal for } {\boldsymbol \mu}\}$.

Conversely, we define a single-valued mapping
\begin{align}
\Lambda : \Pi \to \M \nonumber
\end{align}
as follows: given $\pi \in \Pi$, the state-measure flow $\bmu := \Lambda(\pi)$ is constructed recursively as
\begin{align}
\mu_{t+1}(\,\cdot\,) = \int_{\sX \times \sA} p(\,\cdot\,|x(t),a(t),\mu_t) P^{\pi}(da(t)|x(t)) \mu_t(dx(t)), \nonumber
\end{align}
where $P^{\pi}(da(t)|x(t))$ denotes the conditional distribution of $a(t)$ given $x(t)$ under $\pi$ and $(\mu_{\tau})_{0\leq\tau\leq t}$. Note that if $\pi$ is a Markov policy (i.e., $\pi_t(da(t)|g(t)) = \pi_t(da(t)|x(t))$ for all $t$), then $P^{\pi}(da(t)|x(t)) = \pi_t(da(t)|x(t))$. Using $\Psi$ and $\Lambda$, we now introduce the notion of an equilibrium for the mean-field game.

\begin{definition}
A pair $(\pi,{\boldsymbol \mu}) \in \Pi \times \M$ is a \emph{mean-field equilibrium} if $\pi \in \Phi({\boldsymbol \mu})$ and $\bmu = \Lambda(\pi)$.
\end{definition}

The following structural result shows that the restriction to Markov policies entails no loss of optimality.
\begin{proposition}\label{newprop1}
For any state measure flow $\bmu \in \M$, we have
\begin{align}
\inf_{\pi \in \Pi} J_{\bmu}(\pi) = \inf_{\pi \in \sM} J_{\bmu}(\pi). \nonumber
\end{align}
Furthermore, we have $\Lambda(\Pi) = \Lambda(\sM)$; that is, for any $\pi \in \Pi$, there exists $\hat{\pi} \in \sM$ such that $\mu_t^{\pi} = \mu_t^{\hat{\pi}}$ for all $t\geq0$.
\end{proposition}
\proof{}
This can be proved as in \cite[Proposition 3.2]{SaBaRa17}. Hence, we omit the details.\Halmos
\endproof

This result implies that we can restrict ourselves to Markov policies in the definitions of $\Phi$ and $\Lambda$ without any loss of generality. Furthermore, it also implies that, unlike the finite-population case, the policy in the mean-field equilibrium constitutes a true Nash equilibrium of Markovian type for the infinite-population game problem.

In this section, the main goal is to establish the existence of a mean-field equilibrium. To that end, we impose the assumptions below on the components of the mean-field game model.

\begin{assumption}\label{as1}
\begin{itemize}
\item [ ]
\item [(a)] The cost function $c$ is bounded and continuous with $\| c \| \leq K$.
\item [(b)] The stochastic kernel $p$ is weakly continuous in $(x,a,\mu)$; i.e., $p(\,\cdot\,|x(k),a(k),\mu_k) \rightarrow p(\,\cdot\,|x,a,\mu)$ weakly when $(x(k),a(k),\mu_k) \rightarrow (x,a,\mu)$.
\item [(c)] $\sA$ is compact.
\item [(d)] There exist a constant $\alpha \ge 0$  and a continuous moment function $w: \sX \rightarrow [1,\infty)$  (see \cite[Definition E.7]{HeLa96}) such that
\begin{align}
\sup_{(a,\mu) \in \sA \times \P(\sX)} \int_{\sX} w(y) p(dy|x,a,\mu) \leq \alpha w(x).
\end{align}
\item [(e)] The initial probability measure $\mu_0$ satisfies
\begin{align}
\int_{\sX} w(x) \mu_0(dx) \eqqcolon M < \infty. \nonumber
\end{align}
\end{itemize}
\end{assumption}

\begin{remark}
To ease the exposition of the paper, we assume that the one-stage cost function $c$ is bounded. However, one can extend all the results in this paper to the case where the one-stage cost function $c$ is unbounded with some growth condition (see \cite[Section 2.1]{SaBaRa17}).
\end{remark}

In this section, we prove the existence of a mean-field equilibrium under Assumption~\ref{as1}. Later we will show that this mean-field equilibrium constitutes an approximate Nash-equilibrium for games with sufficiently many agents.

\begin{theorem}\label{thm:MFE} Under Assumption~1, the mean-field game $\bigl( \sX, \sA, p, c, \mu_0 \bigr)$ admits a mean-field equilibrium $(\pi,\bmu)$.
\end{theorem}

The proof of Theorem~\ref{thm:MFE} is given in the next section.

\section{Proof of Theorem~\ref{thm:MFE}}\label{sec3-1}

\subsection{Risk-Sensitive Nonhomogeneous MDPs}\label{risk-MDP}

Note that, given any measure flow $\bmu \in \M$, the optimal control problem for the mean-field game reduces to one of finding an optimal policy for a risk-sensitive nonhomogeneous Markov decision process. Hence, before starting the proof of Theorem~\ref{thm:MFE}, we first derive some relevant results on risk-sensitive \emph{nonhomogeneous} Markov decision processes which are mostly absent and dispersed in the literature. To this end, fix any $\bmu \in \M$ and consider the corresponding optimal control problem. In the remainder of this section, when we say policy, it should be understood as Markov policy as it is sufficient for optimality by Proposition~\ref{newprop1}.

Define
\begin{align}
p_t(\,\cdot\,|x,a) &\coloneqq p(\,\cdot\,|x,a,\mu_t) \nonumber \\
c_t(x,a) &\coloneqq c(x,a,\mu_t). \nonumber
\end{align}
For any $\gamma > 0$, we define
\begin{align}
J_k^n(\pi,x,\gamma) &\coloneqq E^{\pi} \biggl[ e^{\gamma \sum_{t=k}^{n} \beta^{t-k} c_t(x(t),a(t))} \bigg| x(k) = x \biggr] \nonumber \\
J_k^n(x,\gamma) &\coloneqq \inf_{\pi} J_k^n(\pi,x,\gamma) \nonumber \\
J_k(\pi,x,\gamma) &\coloneqq E^{\pi} \biggl[ e^{\gamma \sum_{t=k}^{\infty} \beta^{t-k} c_t(x(t),a(t))} \bigg| x(k) = x \biggr] \nonumber \\
J_k(x,\gamma) &\coloneqq \inf_{\pi} J_k(\pi,x,\gamma), \nonumber
\end{align}
where $n \geq k \geq 0$. With these definitions, the cost function of any policy $\pi$ with initial distribution $\mu_0 = \delta_x$ is $J_0(\pi,x,\lambda)$ and the optimal value function is $J_0(x,\lambda)$. Furthermore, for any $k\geq1$, we have
\begin{align}
E^{\pi} \biggl[ e^{\lambda \sum_{t=k}^{\infty} \beta^{t} c_t(x(t),a(t))} \bigg| x(k) = x \biggr] &= J_k(\pi,x,\lambda \beta^k) \nonumber \\
\inf_{\pi} E^{\pi} \biggl[ e^{\lambda \sum_{t=k}^{\infty} \beta^{t} c_t(x(t),a(t))} \bigg| x(k) = x \biggr] &= J_k(x,\lambda \beta^k). \nonumber
\end{align}

We first prove that finite-horizon cost functions and optimal value functions are continuous and bounded.

\begin{lemma}\label{lemma1}
For each $n \geq k \geq 0$, we have $J_k^n(\cdot,\lambda\beta^k) \in C_b(\sX)$ with $\|J_k^n(\cdot,\lambda\beta^k)\| \leq e^{\lambda K \zeta_{k,n}}$, where $K$ is the constant in Assumption~\ref{as1}-(a) and $\zeta_{k,n} \coloneqq \sum_{t=k}^n \beta^k$ which is less than $\frac{1}{1-\beta}$.
\end{lemma}

\proof{}
Fix any $n$. Then, we prove the result by backward induction on $k$. For $k=n$, we have
\begin{align}
J_n^n(x,\lambda\beta^n) &= \inf_{\pi} E^{\pi} \biggl[ e^{\lambda \beta^n c_n(x(n),a(n))} \bigg| x(n) = x \biggr] \nonumber \\
&= \inf_{a \in \sA} e^{\lambda \beta^n c_n(x,a)}. \nonumber
\end{align}
Since $\sA$ is compact and $e^{\lambda \beta^n c_n(x,a)} \in C_b(\sX \times \sA)$, we have $J_n^n(\cdot,\lambda\beta^n) \in C_b(\sX)$ and obviously $\|J_n^n(\cdot,\lambda\beta^k)\| \leq e^{\lambda K \zeta_{n,n}}$. Suppose the claim holds for $k+1$ and consider $k$. We have
\begin{align}
&J_k^n(x,\lambda\beta^k) = \inf_{\pi} E^{\pi} \biggl[ e^{\lambda\beta^k \sum_{t=k}^{n} \beta^{t-k} c_t(x(t),a(t))} \bigg| x(k) = x \biggr] \nonumber \\
&= \inf_{\pi} E^{\pi} \biggl[ e^{\lambda\beta^k c_k(x(k),a(k))} e^{\lambda\beta^k \sum_{t=k+1}^{n} \beta^{t-k} c_t(x(t),a(t))} \bigg| x(k) = x \biggr] \nonumber \\
&\overset{(a)}{=} \inf_{a \in \sA} \biggl\{  e^{\lambda\beta^k c_k(x,a)}  \inf_{\pi} E^{\pi} \biggl[e^{\lambda\beta^{k+1} \sum_{t=k+1}^{n} \beta^{t-k-1} c_t(x(t),a(t))} \bigg| x(k) = x, a(k) = a \biggr] \biggr\} \nonumber \\
&= \inf_{a \in \sA} \biggl\{ e^{\lambda\beta^k c_k(x,a)} \int_{\sX} J_{k+1}^n(y,\lambda\beta^{k+1}) p_k(dy|x,a) \biggr\}, \nonumber
\end{align}
where (a) follows from the fact that the current state and action pair $(x,a)$ affects the future cost only through the transition probability $p_k(\,\cdot\,|x,a)$ and the multiplicative constant $e^{\lambda\beta^k c_k(x,a)}$. Since $\sA$ is compact, $c_k \in C_b(\sX \times \sA)$, $J_{k+1}^n(\cdot,\lambda\beta^{k+1}) \in C_b(\sX)$, and $p_k$ is weakly continuous, we have $J_{k}^n(\cdot,\lambda\beta^{k}) \in C_b(\sX)$ with $\|J_k^n(\cdot,\lambda\beta^k)\| \leq e^{\lambda K \zeta_{k,n}}$.\Halmos
\endproof

For each $k\geq0$,  define the operator $T_k: C_b(\sX) \rightarrow C_b(\sX)$ as
\begin{align}
[T_ku] (x) \coloneqq \inf_{a \in \sA} \biggl[ e^{\lambda\beta^k c_k(x,a)} \int_{\sX} u(y) p_k(dy|x,a) \biggr] \label{opt:eq}.
\end{align}
Therefore, in view of the proof of Lemma~\ref{lemma1}, for each $n > k \geq 0$,  we have
\begin{align}
\bigl[T_kJ_{k+1}^n(\cdot,\lambda\beta^{k+1})\bigr](x) = J_k^n(x,\lambda\beta^k). \label{opt:eq2}
\end{align}
A similar conclusion can be derived for $n = \infty$ as it is stated in the next lemma.

\begin{lemma}\label{lemma2}
For each $k\geq0$, $\bigl[T_kJ_{k+1}(\cdot,\lambda\beta^{k+1})\bigr](x) = J_k(x,\lambda\beta^k)$.
\end{lemma}

The next lemma states that the finite-horizon optimal cost function converges uniformly to the infinite-horizon optimal cost function. To prove the result, we use the following inequality for $y_1 >0, y_2 >0$:
\begin{align}
e^{\lambda(y_1+y_2)} \leq e^{\lambda y_1} + \lambda e^{\lambda(y_1+y_2)} y_2. \label{ineq}
\end{align}
This inequality is a special case of the inequality $U(y_1+y_2) \leq U(y_1) + U^{'}(y_1+y_2) y_2$ for $U:\R_{+} \rightarrow \R$ convex (see \cite{BaRi14}).

\begin{lemma}\label{lemma3}
For each $k\geq0$, we have
\begin{align}
\bigl\| J_k^n(\cdot,\lambda\beta^k) - J_k(\cdot,\lambda\beta^k) \bigr\| \leq L_k \beta^{n+1} \rightarrow  0 \text{ } \text{ as $n \rightarrow \infty$}, \nonumber
\end{align}
where $L_k \coloneqq \frac{\lambda K}{1-\beta} e^{\frac{\lambda\beta^k K}{1-\beta}}$. Therefore, for each $k \geq 0$, $J_k(\cdot,\lambda\beta^k) \in C_b(\sX)$ with $\|J_k(\cdot,\lambda\beta^k)\| \leq e^{\lambda \frac{K}{1-\beta}}$.
\end{lemma}

\proof{}
Note that $J_k^n(\cdot,\lambda\beta^k) \leq J_k^{n+1}(\cdot,\lambda\beta^k) \leq J_k(\cdot,\lambda\beta^k)$ for all $n$ as $c_t \geq 0$ for all $t\geq0$. Moreover, for any policy $\pi$, we have
\begin{align}
&J_k(\pi,x,\lambda\beta^k) = E^{\pi} \biggl[ e^{\lambda\beta^k \sum_{t=k}^{\infty} \beta^{t-k} c_t(x(t),a(t))} \bigg| x(k) = x \biggr] \nonumber \\
&= E^{\pi} \biggl[ e^{\lambda\beta^k \bigl( \sum_{t=k}^{n} \beta^{t-k} c_t(x(t),a(t)) + \sum_{t=n+1}^{\infty} \beta^{t-k} c_t(x(t),a(t)) \bigr)} \bigg| x(k) = x \biggr] \nonumber \\
&\leq E^{\pi} \biggl[ e^{\lambda\beta^k \sum_{t=k}^{n} \beta^{t-k} c_t(x(t),a(t))} + \lambda\beta^k e^{\lambda\beta^k \sum_{t=k}^{\infty} \beta^{t-k} c_t(x(t),a(t))} \sum_{t=n+1}^{\infty} \beta^{t-k} c_t(x(t),a(t)) \bigg| x(k) = x \biggr] \text{ } \text{(by (\ref{ineq}))}\nonumber \\
&\leq J_k^n(\pi,x,\lambda\beta^k) + \frac{\lambda K}{1-\beta} e^{\frac{\lambda\beta^k K}{1-\beta}} \beta^{n+1} \nonumber \\
&\eqqcolon J_k^n(\pi,x,\lambda\beta^k) + L_k \beta^{n+1}. \nonumber
\end{align}
Therefore, we have
\begin{align}
\bigl\| J_k^n(\cdot,\lambda\beta^k) - J_k(\cdot,\lambda\beta^k) \bigr\| \leq L_k \beta^{n+1} \rightarrow 0 \text{ } \text{ as } n \rightarrow \infty. \nonumber
\end{align}
\Halmos
\endproof

The following is a key result that characterizes the optimal policies. We use this to establish the existence of mean-field equilibrium.

\begin{theorem}\label{opt:thm}
A policy $\pi$ is optimal if and only if, for all $k$,
\begin{align}
&\nu_k^{\pi} \biggl( \biggl\{ (x,a): e^{\lambda\beta^kc_k(x,a)} \int_{\sX} J_{k+1}(y,\lambda\beta^{k+1}) p_k(dy|x,a) \nonumber \\
&\phantom{xxxxxxxxxxxxxxxxxxxxxxxxxxx}= \bigl[ T_k J_{k+1}(\cdot,\lambda\beta^{k+1}) \bigr](x) \biggr\} \biggl) = 1 \label{opt:cond}
\end{align}
where $\nu_k^{\pi} = \mathcal{L}(x(k),a(k))$ under $\pi$ and $\mu_0$.
\end{theorem}

\proof{}
Suppose first that $\pi$ is optimal. Then, one can prove that
\begin{align}
J_k(\cdot,\lambda\beta^k) = J_k(\pi,\cdot,\lambda\beta^k), \text{ } \text{$P^{\pi}$-a.s.}, \label{aux2}
\end{align}
for all $k\geq0$. Moreover, we have
\begin{align}
E^{\pi} \bigl[ J_k(\pi,x(k),\lambda\beta^k) \bigr] &= \int_{\sX \times \sA} e^{\lambda\beta^k c_k(x,a)} \int_{\sX} J_{k+1}(\pi,y,\lambda\beta^{k+1}) p_k(dy|x,a) \nu_k^{\pi}(dx,da) \nonumber \\
&= E^{\pi} \bigl[ J_k(x(k),\lambda\beta^k) \bigr] \text{ } \text{(by (\ref{aux2}))}\nonumber \\
&= \int_{\sX \times \sA} \bigl[T_kJ_{k+1}(\cdot,\lambda\beta^{k+1})\bigr](x) \nu_k^{\pi}(dx,da). \nonumber
\end{align}
Since $e^{\lambda\beta^k c_k(x,a)} \int_{\sX} J_{k+1}(\pi,y,\lambda\beta^{k+1}) p_k(dy|x,a) \geq \bigl[T_kJ_{k+1}(\cdot,\lambda\beta^{k+1})\bigr](x)$ for all $(x,a) \in \sX \times \sA$, (\ref{opt:cond}) holds.

Conversely, suppose that $\pi = \{\pi_0,\pi_1,\ldots\}$ is a policy that satisfies (\ref{opt:cond}). For any $k\geq1$, let $\pi^k = \{\pi_k,\pi_{k+1},\ldots\}$. First, we prove by induction on $n$ that for each $k=0,1,2,\ldots$
\begin{align}
J_k^{n+k}(\pi^k,x,\lambda\beta^k) \leq J_k(x,\lambda\beta^k), \text{ } \text{$P^{\pi}$-a.s.} \label{aux1}
\end{align}
For $n=0$, we have that $J_k^k(\pi^k,x,\lambda\beta^k) = \int_{\sA} e^{\lambda\beta^k c_k(x,a)} \pi_k(da|x)$. On the other hand, we have
\begin{align}
J_k(x,\lambda\beta^k) &= \inf_{a \in \sA} \biggl[ e^{\lambda\beta^k c_k(x,a)} \int_{\sX} J_{k+1}(y,\lambda\beta^{k+1}) p_k(dy|x,a) \biggr] \nonumber \\
&= \int_{\sA} e^{\lambda\beta^k c_k(x,a)} \int_{\sX} J_{k+1}(y,\lambda\beta^{k+1}) p_k(dy|x,a) \pi_k(da|x), \text{ } \text{$P^{\pi}$-a.s.} \nonumber \\
&\geq J_k^k(\pi^k,x,\lambda\beta^k) \text{ } \text{( as $J_{k+1}(y,\lambda\beta^{k+1}) \geq 1$)}. \nonumber
\end{align}
Hence the claim holds for $n=0$. Suppose that the claim holds for $n-1\geq0$; that is, for each $k\geq0$, we have
\begin{align}
J_k^{(n-1)+k}(\pi^k,x,\lambda\beta^k) \leq J_k(x,\lambda\beta^k), \text{ } \text{$P^{\pi}$-a.s.} \nonumber
\end{align}
Then, for each $k\geq0$, we have
\begin{align}
&J_k^{n+k}(\pi^k,x,\lambda\beta^k) = \int_{\sA} e^{\lambda\beta^k c_k(x,a)} \int_{\sX} J_{k+1}^{n+k}(\pi^{k+1},y,\lambda\beta^{k+1}) p_k(dy|x,a) \pi_k(da|x) \nonumber \\
&\leq \int_{\sA} e^{\lambda\beta^k c_k(x,a)} \int_{\sX} J_{k+1}(\pi^{k+1},y,\lambda\beta^{k+1}) p_k(dy|x,a) \pi_k(da|x), \text{ } \text{$P^{\pi}$-a.s.} \text{ } \text{(by induction hypothesis)} \nonumber \\
&= \inf_{a \in \sA} \biggl[ e^{\lambda\beta^k c_k(x,a)} \int_{\sX} J_{k+1}(\pi^{k+1},y,\lambda\beta^{k+1}) p_k(dy|x,a) \biggr], \text{ } \text{$P^{\pi}$-a.s.} \nonumber \\
&\eqqcolon J_k(x,\lambda\beta^k). \nonumber
\end{align}
Hence, this completes the proof of (\ref{aux1}). Note that if $k=0$, then we have $J_0^n(\pi,x,\lambda) \leq J_0(x,\lambda)$ $\mu_0$-a.s. for all $n$. But, by Lemma~\ref{lemma3}, $J_0^n(\pi,x,\lambda) \rightarrow J_0(\pi,x,\lambda)$ uniformly as $n\rightarrow\infty$. Thus, $J_0(\pi,x,\lambda) \leq J_0(x,\lambda)$ $\mu_0$-a.s., and so, $\pi$ is optimal.\Halmos
\endproof

\subsection{Proof of Existence of Mean-Field Equilibrium}\label{existence-mfe}

We are now ready to prove Theorem~\ref{thm:MFE}. For each $t\geq0$, define the following sets
\begin{align}
\P_w^t(\sX) &\coloneqq \{\mu \in \P(\sX): \mu(w) \leq \alpha^t M\} \nonumber \\
\P_w^t(\sX \times \sA) &\coloneqq \{\nu \in \P(\sX \times \sA): \nu_1 \in \P_w^t(\sX)\}, \nonumber
\end{align}
where for any $\nu \in \P(\sX \times \sA)$, $\nu_1$ denotes the marginal of $\nu$ on $\sX$. Let $\Xi \coloneqq \prod_{t=0}^{\infty} \P_w^t(\sX \times \sA)$, which is equipped with the infinite product topology generated by weak convergence.

For any $\bnu \in \Xi$ and $t\geq0$, we define the operator $T_t^{\bnu}: C_b(\sX) \rightarrow C_b(\sX)$ as
\begin{align}
[T_t^{\bnu}u](x) = \inf_{a \in \sA} \biggl[ e^{\lambda \beta^t c(x,a,\nu_{t,1})} \int_{\sX} u(y) p(dy|x,a,\nu_{t,1}) \biggr]. \nonumber
\end{align}
We let $J_{*,t}^{\bnu}(\cdot,\lambda\beta^t): \sX \rightarrow \R$ denote the risk-sensitive optimal value function at time $t$ of the nonhomogeneous Markov decision process with the one-stage cost functions $\bigl\{c(x,a,\nu_{t,1})\bigr\}_{t\geq0}$ and the transition probabilities $\bigl\{p(\,\cdot\,|x,a,\nu_{t,1})\bigr\}_{t\geq0}$. By Lemma~\ref{lemma3}, for all $t\geq0$, $J_{*,t}^{\bnu}(\cdot,\lambda\beta^t)$ is continuous and bounded.

To prove the existence of a mean-field equilibrium, we use the technique in our previous paper \cite{SaBaRa17} adopted from \cite{JoRo88}, which enables us to transform the fixed point equation $\pi \in \Psi(\Lambda(\pi))$ characterizing the mean-field equilibrium into a fixed point equation of a set-valued mapping from the set of state-action measure flows $\P(\sZ \times \sA)^{\infty}$ into itself. Then, using Kakutani's fixed point theorem (\cite[Corollary 17.55]{AlBo06}), we deduce the existence of a mean-field equilibrium. We note that the technique used here to prove the existence of a mean-field equilibrium is similar to the one in our previous paper \cite{SaBaRa17}, in which we have studied a risk-neutral version of the same problem. However, there is a crucial difference in the details of the proofs between this problem and the risk-neutral one. In the risk-neutral case, the cost function is in an additive form, and this results in a contractive dynamic programming operator. However, in the risk-sensitive case, the cost function is
in a multiplicative form, and therefore, the dynamic programming operator $T_t$ is not contractive, which complicates the analysis.

As we mentioned above, we first transform the fixed point equation $\pi \in \Psi(\Lambda(\pi))$  into a fixed point equation of a set-valued mapping from $\P(\sX \times \sA)^{\infty}$ into itself. To that end, define the set-valued mapping $\Gamma: \Xi \rightarrow 2^{\P(\sX \times \sA)^{\infty}}$ as follows:
\begin{align}
\Gamma(\bnu) = C(\bnu) \cap B(\bnu), \nonumber
\end{align}
where
\begin{align}
C(\bnu) &\coloneqq \biggl\{ \bnu' \in \P(\sX \times \sA)^{\infty}: \nu'_{0,1} = \mu_0 \text{ and } \nonumber \\
&\phantom{xxxxxxxxxxxxxx}\nu'_{t+1,1}(\,\cdot\,) = \int_{\sX \times \sA} p(\,\cdot\,|x,a,\nu_{t,1}) \nu_t(dx,da) \biggr\} \nonumber \\
\intertext{and}
B(\bnu) &\coloneqq \biggl\{ \bnu' \in \P(\sX \times \sA)^{\infty}: \forall t\geq0, \text{ } \nu_t' \biggl( \biggr\{ (x,a) : e^{\lambda\beta^t c(x,a,\nu_{t,1})} \nonumber \\
&\phantom{xxxx} \int_{\sX} J_{*,t+1}^{\bnu}(y,\lambda\beta^{t+1}) p(dy|x,a,\nu_{t,1}) = T_t^{\bnu} J^{\bnu}_{*,t+1}(x,\lambda\beta^{t+1}) \biggr\} \biggr) = 1 \biggr\}. \nonumber
\end{align}
Note that the set $C(\bnu)$ characterizes the consistency of the mean-field term with the distribution of a generic agent, and the set $B(\bnu)$ characterizes optimality of the policy that is obtained by disintegrating the state-action measure-flow, for the mean-field term. The following proposition implies that the image of $\Xi$ under $\Gamma$ is contained in $2^{\Xi}$.

\begin{proposition}({\cite[Proposition 3.7]{SaBaRa17}})\label{prop2}
For any $\bnu \in \Xi$, we have $\Gamma(\bnu) \subset \Xi$.
\end{proposition}


We say that $\nu \in \Xi$ is a fixed point of $\Gamma$ if $\bnu \in \Gamma(\bnu)$. The following proposition makes the connection between mean-field equilibria and the fixed points of $\Gamma$.

\begin{proposition}\label{prop1}
Suppose that $\Gamma$ has a fixed point $\bnu = (\nu_t)_{t \ge 0}$. Construct a Markov policy $\pi = (\pi_t)_{t \ge 0}$ by disintegrating each $\nu_t$ as $\nu_t(dx,da) = \nu_{t,1}(dx) \pi_t(da|x)$, and let $\bnu_1 = (\nu_{t,1})_{t \ge 0}$. Then the pair $(\pi,\bnu_1)$ is a mean-field equilibrium.
\end{proposition}

\proof{}
If $\bnu \in \Gamma(\bnu)$, then the corresponding Markov policy $\pi$ satisfies (\ref{opt:cond}) for $\bnu$. Therefore, by Theorem~\ref{opt:thm}, $\pi \in \Phi(\bnu_1)$. Moreover, since $\bnu \in C(\bnu)$, we have $\Lambda(\pi) = \nu_1$. This completes the proof.\Halmos
\endproof

\noindent By Proposition~\ref{prop1}, it suffices to prove that $\Gamma$ has a fixed point in order to establish the existence of a mean-field equilibrium. To prove this, we use Kakutani's fixed point theorem (\cite[Corollary 17.55]{AlBo06}). Note that, since $w$ is a continuous moment function, the set $\P^t_w(\sX)$ is compact with respect to the weak topology (\cite[Proposition E.8, p. 187]{HeLa96}), and so, $\P^t_w(\sX \times \sA)$ is tight as $\sA$ is compact. Furthermore, $\P^t_w(\sX \times \sA)$ is closed with respect to the weak topology. Hence, $\P^t_w(\sX \times \sA)$ is compact with respect to the weak topology. Therefore, $\Xi$ is compact with respect to the product topology. In addition, $\Xi$ is also convex. Furthermore, it can be proved that $C(\bnu) \cap B(\bnu) \neq \emptyset$ for any $\bnu \in \Xi$. We can also show that both $C(\bnu)$ and $B(\bnu)$ are convex, and thus their intersection is also convex. Moreover, $\Xi$ is a convex compact subset of a locally convex topological space $\M(\sX \times \sA)^{\infty}$, where $\M(\sX \times \sA)$ denotes the set of finite signed measures on $\sX \times \sA$. The final result we need to prove in order to deduce the existence of a fixed point of $\Gamma$ by an appeal to Kakutani's fixed point theorem is the following:

\begin{proposition}\label{prop3}
The graph of $\Gamma$, i.e., the set
	$$
	\Gr(\Gamma) := \left\{ (\bnu,\bxi) \in \Xi \times \Xi : \bxi \in \Gamma(\bnu)\right\},
	$$
is closed.
\end{proposition}

\proof{}
Let $\bigl\{(\bnu^{(n)},\bxi^{(n)})\bigr\}_{n\geq1} \subset \Xi \times \Xi$ be such that $\bxi^{(n)} \in \Gamma(\bnu^{(n)})$ for all $n$ and $(\bnu^{(n)},\bxi^{(n)}) \rightarrow (\bnu,\bxi)$ as $n\rightarrow\infty$ for some $(\bnu,\bxi) \in \Xi \times \Xi$. To prove $\Gr(\Gamma)$ is closed, it is sufficient to prove $\bxi \in \Gamma(\bnu)$.


One can prove that $\bxi \in C(\bnu)$ (see the proof of \cite[Proposition 3.9]{SaBaRa17}). Hence, it remains to prove that $\bxi \in B(\bnu)$. For each $n$ and $t$, let us define
\begin{align}
F^{(n)}_t(x,a) &= e^{\lambda\beta^t c(x,a,\nu^{(n)}_{t,1})} \int_{\sX} J^{\bnu^{(n)}}_{*,t+1}(y,\lambda\beta^{t+1}) p(dy|x,a,\nu^{(n)}_{t,1}) \nonumber \\
\intertext{and}
F_t(x,a) &= e^{\lambda\beta^t c(x,a,\nu_{t,1})} \int_{\sX} J^{\bnu}_{*,t+1}(y,\lambda\beta^{t+1}) p(dy|x,a,\nu_{t,1}). \nonumber
\end{align}
By definition,
\begin{align}
J^{\bnu^{(n)}}_{*,t}(x,\lambda\beta^t) = \inf_{a \in \sA} F^{(n)}_t(x,a) \text{ } \text{ and } \text{ } J^{\bnu}_{*,t}(x,\lambda\beta^t) = \inf_{a \in \sA} F_t(x,a). \nonumber
\end{align}
By assumption, we have
\begin{align}
1 = \xi^{(n)}_t\biggl( \biggl\{ (x,a): F^{(n)}_t(x,a) = J^{\bnu^{(n)}}_{*,t}(x,\lambda\beta^t) \biggr\} \biggr), \text{ } \text{for all $n$}. \nonumber
\end{align}
Let $A_t^{(n)} \coloneqq \bigl\{ (x,a): F^{(n)}_t(x,a) = J^{\bnu^{(n)}}_{*,t}(x,\lambda\beta^t) \bigr\}$. Since both $F^{(n)}_t$ and $J^{\bnu^{(n)}}_{*,t}$ are continuous, $A_t^{(n)}$ is closed. Define $A_t \coloneqq \bigl\{ (x,a): F_t(x,a) = J^{\bnu}_{*,t}(x,\lambda\beta^t) \bigr\}$ which is also closed as both $F_t$ and $J^{\bnu}_{*,t}$ are continuous.

Suppose that $F_t^{(n)}$ converges to $F_t$ continuously and $J^{\bnu^{(n)}}_{*,t}$ converges to $J^{\bnu}_{*,t}$ continuously, as $n\rightarrow\infty$\footnote{Suppose $g$, $g_n$ ($n\geq1$) are measurable functions on metric space $\sE$. The sequence $g_n$ is said to converge to $g$ continuously if $\lim_{n\rightarrow\infty}g_n(e_n)=g(e)$ for any $e_n\rightarrow e$ where $e \in \sE$.} (see \cite[p. 388]{Ser82}).

For each $M\geq1$, define $B_t^M \coloneqq \bigl\{ (x,a): F_t(x,a) \geq J^{\bnu}_{*,t}(x,\lambda\beta^t) + \epsilon(M) \bigr\}$ which is closed, where $\epsilon(M) \rightarrow 0$ as $M\rightarrow \infty$. Since both $F_t$ and $J^{\bnu}_{*,t}$ are continuous, we can choose $\{\epsilon(M)\}_{M\geq1}$ so that $\xi_t(\partial B_t^M) = 0$ for each $M$. Since $A_t^c = \bigcup_{M=1}^{\infty} B_t^M$ and $B_t^M \subset B_t^{M+1}$, we have by the monotone convergence theorem
\begin{align}
\xi^{(n)}_t\big(A_t^c \cap A_t^{(n)}\big) = \liminf_{M\to\infty} \xi^{(n)}_t\big(B^M_t \cap A_t^{(n)}). \nonumber
\end{align}
Hence, we have
\begin{align}
1 &= \limsup_{n\rightarrow\infty} \liminf_{M\rightarrow\infty} \biggl\{ \xi^{(n)}_t\big(A_t \cap A^{(n)}_t\big) + \xi^{(n)}_t\big(B^M_t \cap A_t^{(n)}\big)\biggr\} \nonumber\\
&\leq \liminf_{M\rightarrow\infty} \limsup_{n\rightarrow\infty}  \biggl\{\xi^{(n)}_t\big(A_t \cap A^{(n)}_t\big) + \xi^{(n)}_t\big(B^M_t \cap A_t^{(n)}\big)\biggr\}. \nonumber
\end{align}
\noindent For fixed $M$, we prove that the second term in the last expression converges to zero as $n\rightarrow\infty$. First, note that $\xi^{(n)}_t$ converges weakly to $\xi_t$ as $n\rightarrow\infty$ when both measures are restricted to $B_t^M$, as $B_t^M$ is closed and $\xi_t(\partial B_t^M)=0$ \cite[Theorem 8.2.3]{Bog07}. Furthermore, $1_{A^{(n)}_t \cap B^M_t}$ converges continuously to $0$: if $(x^{(n)},a^{(n)}) \rightarrow (x,a)$ in $B_t^M$, then
\begin{align}
\lim_{n\rightarrow\infty} F_t^{(n)}(x^{(n)},a^{(n)}) &= F_t(x,a) \nonumber \\
&\geq J_{*,t}(x,\lambda\beta^t) + \epsilon(M) \nonumber \\
&= \lim_{n\rightarrow\infty} J^{(n)}_{*,t}(x^{(n)},\lambda\beta^t) + \epsilon(M). \nonumber
\end{align}
Hence, for large enough $n$, we have $F_t^{(n)}(x^{(n)},a^{(n)}) > J^{(n)}_{*,t}(x^{(n)},\lambda\beta^t)$, which implies that $(x^{(n)},a^{(n)}) \not\in A^{(n)}_t$. Then, by \cite[Theorem 3.5]{Lan81}, for each $M$ we have
\begin{align*}
\limsup_{n\rightarrow\infty} \xi^{(n)}_t\big(B^M_t \cap A^{(n)}_t\big) = 0.
\end{align*}
Therefore, we obtain
\begin{align*}
1 &\leq  \limsup_{n\rightarrow\infty} \xi^{(n)}_t\big(A_t \cap A_t^{(n)}\big)\\
&\le \limsup_{n\rightarrow\infty} \xi_t^{(n)}(A_t) \nonumber \\
&\leq \xi_t(A_t), \nonumber
\end{align*}
where the last inequality follows from the portmanteau theorem \cite[Theorem 2.1]{Bil99} and the fact that $A_t$ is closed. Hence, $\xi_t(A_t)=1$. Since $t$ is arbitrary, this is true for all $t$. This means that $\bxi \in B(\bnu)$. Therefore, $\bxi \in \Gamma(\bnu)$ which completes the proof under the assumption that $F_t^{(n)}$ converges to $F_t$ continuously and $J^{\bnu^{(n)}}_{*,t}$ converges to $J^{\bnu}_{*,t}$ continuously, as $n\rightarrow\infty$, which we prove next.\Halmos
\endproof

Note that, for continuous functions, continuous convergence coincides with uniform convergence over compact sets (see \cite[Lemma 2.1]{Lan81}). Therefore, it suffices to establish that $F^{(n)}_t$ uniformly converges to $F_t$ over compact sets and $J^{\bnu^{(n)}}_{*,t}$ uniformly converges to $J^{\bnu}_{*,t}$ over compact sets, as these functions are all continuous.
Furthermore, if $J^{\bnu^{(n)}}_{*,t}$ converges to $J^{\bnu}_{*,t}$ continuously for all $t$, then $F^{(n)}_t$ also converges to $F_t$ continuously for all $t$. Indeed, let $(x^{(n)},a^{(n)}) \rightarrow (x,a)$. Since $J^{\bnu^{(n)}}_{*,t+1}(y,\lambda\beta^{t+1}) \leq e^{\lambda\frac{K}{1-\beta}}$ for all $n\geq1$, by \cite[Theorem 3.5]{Lan81} we have
\begin{align}
\lim_{n\rightarrow\infty} F_t^{(n)}(x^{(n)},a^{(n)}) &= \lim_{n\rightarrow\infty} \biggl[ e^{\lambda\beta^t c(x^{(n)},a^{(n)},\nu^{(n)}_{t,1})} \int_{\sX} J^{\bnu^{(n)}}_{*,t+1}(y,\lambda\beta^{t+1}) p(dy|x^{(n)},a^{(n)},\nu^{(n)}_{t,1}) \biggr] \nonumber \\
&= e^{\lambda\beta^t c(x,a,\nu_{t,1})} \int_{\sX} J^{\bnu}_{*,t+1}(y,\lambda\beta^{t+1}) p(dy|x,a,\nu_{t,1}) \nonumber \\
&= F_t(x,a). \nonumber
\end{align}
Therefore, it is sufficient to prove that $J^{\bnu^{(n)}}_{*,t}$ uniformly converges to $J^{\bnu}_{*,t}$ over compact sets, for all $t$.

\begin{proposition}\label{prop4}
For any compact $K\subset \sX$, we have
\begin{align}
\lim_{n\rightarrow\infty} \sup_{x \in K} \bigl| J^{\bnu^{(n)}}_{*,t}(x,\lambda\beta^t) - J^{\bnu}_{*,t}(x,\lambda\beta^t)| = 0 \nonumber
\end{align}
for all $t\geq0$. Therefore, $J^{\bnu^{(n)}}_{*,t}$ converges to $J^{\bnu}_{*,t}$ continuously as $n\rightarrow\infty$, for all $t$.
\end{proposition}

\proof{}

For each $k \geq t \geq0$, $\bnu \in \Xi$, we let $J_{*,t}^{\bnu}(x,\lambda\beta^t,k)$ denote the optimal value function at time $t$ of the risk-sensitive nonhomogeneous MDP with finite horizon $k$ corresponding to $\bnu$. By Lemma~\ref{lemma3}, we have
\begin{align}
\| J_{*,t}^{\bnu}(\cdot,\lambda\beta^t,k) - J_{*,t}^{\bnu}(\cdot,\lambda\beta^t) \| \leq L_t \beta^{k+1}. \label{aux3}
\end{align}
We first prove the following result for finite horizon optimal value functions.

\begin{lemma}\label{lemma4}
For all $t\geq0$ and for any compact $K \subset \sX$, we have
\begin{align}
\lim_{n\rightarrow\infty} \sup_{x \in K} \bigl| J_{*,t}^{\bnu^{(n)}}(x,\lambda\beta^t,m+t) -  J_{*,t}^{\bnu}(x,\lambda\beta^t,m+t) \bigr| = 0 \text{  } \text{ for any $m\geq0$}. \nonumber
\end{align}
\end{lemma}

\proof{}
We prove the result by induction on $m$. For $m=0$, we have
\begin{align}
\sup_{x \in K} \bigl| J_{*,t}^{\bnu^{(n)}}(x,\lambda\beta^t,t) -  J_{*,t}^{\bnu}(x,\lambda\beta^t,t) \bigr|  &= \sup_{x \in K} \biggl| \inf_{a \in \sA} e^{\lambda\beta^t c(x,a,\nu^{(n)}_{t,1})} -  \inf_{a \in \sA} e^{\lambda\beta^t c(x,a,\nu_{t,1})}\biggr| \nonumber \\
&\leq \sup_{(x,a) \in K \times \sA} \biggl| e^{\lambda\beta^t c(x,a,\nu^{(n)}_{t,1})} -  e^{\lambda\beta^t c(x,a,\nu_{t,1})}\biggr|. \nonumber
\end{align}
Since $e^{\lambda\beta^t c(\cdot,\cdot,\nu^{(n)}_{t,1})}$ converges to $e^{\lambda\beta^t c(\cdot,\cdot,\nu_{t,1})}$ continuously and since $e^{\lambda\beta^t c(\cdot,\cdot,\nu_{t,1})}$ is continuous, the last term converges to $0$ as continuous convergence and uniform convergence on compact sets are equivalent in this case. Hence, the claim holds for $m=0$.

Suppose the claim holds for a general $m$. Then consider $m+1$. We have
\begin{align}
&\sup_{x \in K} \bigl| J_{*,t}^{\bnu^{(n)}}(x,\lambda\beta^t,m+1+t) -  J_{*,t}^{\bnu}(x,\lambda\beta^t,m+1+t) \bigr| \nonumber \\
&= \sup_{x \in K} \bigg| \inf_{a \in \sA} \biggl[ e^{\lambda\beta^t c(x,a,\nu^{(n)}_{t,1})} \int_{\sX} J_{*,t+1}^{\bnu^{(n)}}(y,\lambda\beta^{t+1},m+1+t) p(dy|x,a,\nu_{t,1}^{(n)}) \biggr] \nonumber \\
&\phantom{xxxxxxx} -\inf_{a \in \sA} \biggl[ e^{\lambda\beta^t c(x,a,\nu_{t,1})} \int_{\sX} J_{*,t+1}^{\bnu}(y,\lambda\beta^{t+1},m+1+t) p(dy|x,a,\nu_{t,1}) \biggr] \bigg| \nonumber \\
&\leq \sup_{(x,a) \in K \times \sA} \bigg| e^{\lambda\beta^t c(x,a,\nu^{(n)}_{t,1})} \int_{\sX} J_{*,t+1}^{\bnu^{(n)}}(y,\lambda\beta^{t+1},m+1+t) p(dy|x,a,\nu_{t,1}^{(n)}) \nonumber \\
&\phantom{xxxxxxx} - e^{\lambda\beta^t c(x,a,\nu_{t,1})} \int_{\sX} J_{*,t+1}^{\bnu}(y,\lambda\beta^{t+1},m+1+t) p(dy|x,a,\nu_{t,1}) \bigg|. \nonumber
\end{align}
Note that by induction hypothesis and by Lemma~\ref{lemma1}, $J_{*,t+1}^{\bnu^{(n)}}(\cdot,\lambda\beta^{t+1},m+(t+1))$ converges to $J_{*,t+1}^{\bnu}(\cdot,\lambda\beta^{t+1},m+(t+1))$ continuously. Then, by \cite[Theorem 3.5]{Lan81}, the last term converges to zero as continuous convergence is equivalent to uniform convergence on compact sets for continuous functions.\Halmos
\endproof

Using the above lemma, we now complete the proof of Proposition~\ref{prop4}. Fix any compact $K \subset \sX$. For all $m\geq0$, we have
\begin{align}
&\sup_{x \in K} \bigl| J^{\bnu^{(n)}}_{*,t}(x,\lambda\beta^t) - J^{\bnu}_{*,t}(x,\lambda\beta^t)| \nonumber \\
&\leq \| J^{\bnu^{(n)}}_{*,t}(\cdot,\lambda\beta^t) - J^{\bnu^{(n)}}_{*,t}(\cdot,\lambda\beta^t,m+t)\| \nonumber \\
&\phantom{xxxxxxxx}+ \sup_{x \in K} \bigl| J^{\bnu^{(n)}}_{*,t}(\cdot,\lambda\beta^t,m+t) - J^{\bnu}_{*,t}(\cdot,\lambda\beta^t,m+t) \bigr| \nonumber \\
&\phantom{xxxxxxxxxxxxxxxxxxxxxx}+ \|J^{\bnu}_{*,t}(\cdot,\lambda\beta^t) - J^{\bnu}_{*,t}(\cdot,\lambda\beta^t,m+t)\| \nonumber \\
&\leq 2 L_t \beta^{m+t+1} + \sup_{x \in K} \bigl| J^{\bnu^{(n)}}_{*,t}(\cdot,\lambda\beta^t,m+t) - J^{\bnu}_{*,t}(\cdot,\lambda\beta^t,m+t) \bigr| \nonumber
\end{align}
This last expression can be made arbitrarily small by first choosing large enough $m$ and then large enough $n$.\Halmos
\endproof

Recall that $\Xi$ is a compact convex subset of the locally convex topological space $\M(\sX \times \sA)^{\infty}$. Furthermore, $\Gamma$ has closed graph by Proposition~\ref{prop3}, and it takes nonempty convex values. Therefore, by Kakutani's fixed point theorem \cite[Corollary 17.55]{AlBo06}, $\Gamma$ has a fixed point. Then, Theorem~\ref{thm:MFE} follows from Proposition~\ref{prop1}.

\begin{remark}
Note that, given any state measure flow $\bmu \in \M$, it is possible to formulate the corresponding risk-sensitive optimal control problem, defined in Section~\ref{risk-MDP},
as a zero-sum dynamic game with risk-neutral cost function using the following duality relation between logarithmic moment generating function and relative entropy \cite{Jas07}:
\begin{align}
\log\biggl(\int_{\sX} e^{g(x)} \zeta(dx)\biggl) = \sup_{\nu \in \P(\sX)} \biggl[ \int_{\sX} g(x) \nu(dx) - h(\nu \| \zeta) \biggr], \label{ent}
\end{align}
where $g:\sX \rightarrow \R$ is bounded and measurable and $h(\nu \| \zeta)$ denotes the relative entropy of $\nu$ with respect to $\zeta$ \cite{PrMeRu96}.
Indeed, one can prove along the lines of \cite{HeMa96,Jas07} that the risk-sensitive optimal control problem is equivalent to the following risk-neutral zero-sum dynamic game. In this game model, the state space is $\sX$, the action space for Player~1 (minimizer) is $\sA$, and the action space for Player~2 (maximizer) is $\P(\sX)$. For any $(x,a,q) \in \sX \times \sA \times \P(\sX)$, the one-stage cost (reward for Player~2) function at time $t$ is
\begin{align}
d_t(x,a,q) \coloneqq \lambda \beta^t c_t(x,a) - h(q \| p_t(\,\cdot\,|x,a)). \nonumber
\end{align}
The evolution of the system is as follows. At a state $x_t$ ($t\geq0$), Player~1 chooses its action $a_t$ from $\sA$ (possibly randomized) based on the current time $t$ and the current state $x_t$ and Player~2 chooses its action from $\P(\sX)$ based on the current time $t$ and the current state-action pair $(x_t,a_t)$. Let $\pi = \{\pi_t\}$ denote the policy for Player~1 and let $\eta = \{\eta_t\}$ denote the policy for Player~2. Hence, at each time $t$, $\pi_t(da|x)$ is a stochastic kernel on $\sA$ given $\sX$ and $\eta_t(dy|x,a)$ is a stochastic kernel on $\sX$ given $\sX \times \sA$. Then, the system moves to the next state according to the probability distribution $\eta_t(dy|x,a)$; that is, the state transition probability of the game model is given by $r(dy|x,a,q) \coloneqq \delta_{q}(dy)$.

For any policy pair $(\pi,\eta)$ and initial state $x$, the risk-neutral cost of the game is given by
\begin{align}
W_0(x,\pi,\eta) &\coloneqq E^{\pi,\eta} \biggl[ \sum_{t=0}^{\infty}  d_t(x_t,a_t,q_t) \biggr]  \nonumber \\
&= E^{\pi,\eta} \biggl[ \sum_{t=0}^{\infty} \lambda \beta^t c_t(x_t,a_t) - h(q_t(\,\cdot\,|x_t,a_t) \| p_t(\,\cdot\,|x_t,a_t)) \biggr]. \nonumber
\end{align}
Then, the optimal upper-value function is defined as
\begin{align}
W_0(x) \coloneqq \inf_{\pi} \sup_{\eta} W_0(x,\pi,\eta). \nonumber
\end{align}
For $k\geq1$, let $W_k$ denote the optimal upper-value function of the game at time $k$. Then, these functions satisfy the following sequence of Isaacs equations (see \cite[Lemma 3.5]{HeMa96} and \cite[Lemma 1]{Jas07}):
\begin{align}
W_k(x) &= \inf_{a \in \sA} \sup_{q \in \P(\sX)} \biggl[ \int_{\sX} W_{k+1}(y) q(dy) + \lambda \beta^k c_k(x,a) - h(q\| p_k(\,\cdot\,|x,a)) \biggr].  \label{ent-opt}
\end{align}
Note that, by first using duality formula (\ref{ent}) and then taking the exponentials of both sides, we can write (\ref{ent-opt}) as
\begin{align}
e^{W_k(x)} = \inf_{a \in \sX} \biggl[e^{\lambda \beta^k c_k(x,a)} \int_{\sX} e^{W_{k+1}(y)} p_k(dy|x,a)\biggr], \nonumber
\end{align}
which is exactly the dynamic programming recursion for the risk-sensitive optimal control problem with $e^{W_k(x)} = J_k(x,\lambda\beta^k)$ in Lemma~\ref{lemma2}. Then, by using similar arguments as in \cite[Proposition 1 and Lemma 1]{Jas07}, it is possible to deduce that this risk-neutral game model is equivalent to the risk-sensitive optimal control problem in Section~\ref{risk-MDP}.

Equivalency of the two problems suggests that the existence of mean-field equilibrium can be established by analyzing the equivalent zero-sum dynamic game with risk-neutral cost function. To do this, one has to identify the relation between the state-measure flows in zero-sum dynamic game and the state-measure flows in the control problem along with the fact that a state-measure flow in any mean-field equilibrium for risk-neutral zero-sum game corresponds to a state-measure flow in some mean-field equilibrium for the risk-sensitive control problem and vice versa. This is an interesting future research direction.

We note that the same equivalency result may not hold for the finite-population case; that is, one cannot write the problem of a generic agent in the finite-population regime as a zero-sum dynamic game with risk-neutral cost function. One of the reasons for this is that, in the finite-population case, there is no dynamic programming recursion for the optimal value functions because of the coupling between agents and the information structure. However, it is still an interesting research direction to study the finite-population mean-field game in which a generic agent is faced with zero-sum dynamic game of the form described above.

\end{remark}

\section{Existence of Approximate Markov-Nash Equilibria}\label{sec4}

In this section, we prove that the policy generated by the mean-field equilibrium, when adopted by each agent, is approximate Markov-Nash equilibrium for games with a sufficiently large number of agents. Let $(\pi,\bmu)$ denote the pair in the mean-field equilibrium. In addition to Assumption~\ref{as1}, we impose an additional assumption, which is stated below. To this end, let $d_{BL}$ denote the bounded Lipschitz metric (see, e.g., \cite[Proposition 11.3.2]{Dud89}) on $\P(\sX)$ that metrizes the weak topology, and define the following moduli of continuity:
\begin{align}
\omega_{p}(r) &\coloneqq \sup_{(x,a) \in \sX\times\sA} \sup_{\substack{\mu,\nu: \\ d_{BL}(\mu,\nu)\leq r}} \|p(\,\cdot\,|x,a,\mu) - p(\,\cdot\,|x,a,\nu)\|_{TV} \nonumber \\
\omega_{c}(r) &\coloneqq \sup_{(x,a) \in \sX\times\sA} \sup_{\substack{\mu,\nu: \\ d_{BL}(\mu,\nu)\leq r}} |c(x,a,\mu) - c(x,a,\nu)|. \nonumber
\end{align}

\begin{assumption}\label{as2}
\begin{itemize}
\item [ ]
\item [(a)] $\omega_p(r) \rightarrow 0$ and $\omega_c(r) \rightarrow 0$ as $r\rightarrow0$.
\item [(b)] For each $t\geq0$, $\pi_t: \sX \rightarrow \P(\sA)$ is weakly continuous.
\end{itemize}
\end{assumption}

\begin{remark}
Note that, if the state transition probability $p$ is independent of the mean-field term, then Assumption~\ref{as2}-(a) for $p$ is always true.
\end{remark}

\begin{remark}\label{continuity}

To ensure Assumption~\ref{as2}-(b), we can impose the following uniqueness condition, which is common in the mean-field literature (see, e.g., \cite[Assuption 4]{GoMoSo10}, \cite[Assumption A5]{SeCa16}, \cite[Assumption H5]{HuMaCa06}, \cite[Assumption A9]{SeCa16-3}): Suppose that, for each $\bmu$, there exists a unique minimizer $a_x \in \sA$ of
\begin{align}
e^{\lambda \beta^t c(x,\,\cdot\,,\mu_t)} \int_{\sX} J_{*,t+1}^{\bmu}(x') p(dx'|x,\,\cdot\,,\mu_t) \eqqcolon H_t(x,\,\cdot\,), \label{unique}
\end{align}
for each $x \in \sX$ and for all $t\geq0$. Under this condition, one can prove that the policy $\pi$ in mean-field equilibrium is deterministic and weakly continuous. To prove this, fix any $t\geq0$ and consider the policy $\pi_t$ at time $t$ in $\pi$. By (\ref{unique}), we must have $\pi_t(\,\cdot\,|x) = \delta_{f_t(x)}(\,\cdot\,)$ for some $f_t: \sX\rightarrow \sA$ that minimizes $H_{t}(x,\,\cdot\,)$ of the above form. We now prove that $f_t$ is continuous, which implies weak continuity. Suppose $x_n \rightarrow x$ in $\sX$. Note that $g_t(\,\cdot\,) = \min_{a \in \sA} H_{t}(\,\cdot\,,a)$ is continuous since $\sA$ is compact and $H_t$ is continuous. Therefore, every accumulation point of the sequence $\{f_t(x_n)\}_{n\geq1}$ must be a minimizer for $H_{t}(x,\,\cdot\,)$. Since there exists a unique minimizer $f_t(x)$ of $H_{t}(x,\,\cdot\,)$, the set of all accumulation points of $\{f_t(x_n)\}_{n\geq1}$ must be singleton $\{f_t(x)\}$. This implies that $f_t(x_n)$ converges to $f_t(x)$ since $\sA$ is compact. Hence, $f_t$ is continuous. Thus, Assumption~\ref{as2}-(b).

One can establish the existence of a unique minimizer to (\ref{unique}) under the following conditions on the system components. Suppose that $\sX = \R^d$ and $\sA \subset \R^m$ is convex. In addition, suppose that $p(dx'|x,a,\mu) = \varrho(x'|x,a,\mu) m(dx')$, where $m$ denotes the Lebesgue measure. Assume that both $\varrho$ and $c$ are convex in $a$. Hence, for any $a \in \sA$, (\ref{unique}) can be written as
\begin{align}
e^{\lambda \beta^t c(x,\,\cdot\,,\mu_t)}  \int_{\sX} J_{*,t+1}^{\bmu}(x') \varrho(x'|x,\,\cdot\,,\mu_t) m(dx'). \nonumber
\end{align}
Since $c$ and $\varrho$ are convex in $a$, the last expression above is strictly convex in $a$ as exponential function is strictly convex. Hence, there exists a unique minimizer $a_x \in \sA$ for (\ref{unique}).
\end{remark}

The following theorem is the main result of this section:
\begin{theorem}\label{appr-thm}
For any $\varepsilon>0$, there exists $N(\varepsilon)$ such that for $N\geq N(\varepsilon)$, the policy ${\boldsymbol \pi}^{(N)} = (\pi,\ldots,\pi)$ is an $\varepsilon$-Markov-Nash equilibrium for the game with $N$ agents.
\end{theorem}

\subsection*{Proof of Theorem~\ref{appr-thm}}\label{sec4-1}

To prove Theorem~\ref{appr-thm}, we first consider the game problem with a finite-horizon discounted risk-sensitive cost function. For any finite-horizon $n$, we prove the following result.

\begin{theorem}\label{app-mainthm}
For each $N$, let $\tpi^{(N)}$ be some arbitrary weakly continuous Markov policy for Agent~1. Then, the following are true:
\begin{align}
\lim_{N\rightarrow\infty} J_1^{(N),n}(\pi,\ldots,\pi) &= J_{\bmu}^n(\pi), \label{app-main1} \\
\lim_{N\rightarrow\infty} \bigl| J_1^{(N),n}(\tpi^{(N)},\pi,\ldots,\pi)- J_{\bmu}^n(\tpi^{(N)}) \bigr| &= 0. \label{app-main2}
\end{align}
\end{theorem}

In Theorem~\ref{app-mainthm}, superscript $n$ denotes the finite-horizon; that is, for instance, we have
\begin{align}
J_1^{(N),n}(\pi,\ldots,\pi) \coloneqq E^{\bpi^{(N)}} \biggl[ e^{\lambda \sum_{t=0}^{n} \beta^t c(x_1^N(t),a_1^N(t),e_t^{(N)})} \biggr]. \nonumber
\end{align}

Note that, in risk-sensitive cost, the one-stage cost functions are in a multiplicative form as opposed to risk-neutral case, and this makes the analysis of the approximation problem quite complicated. Therefore, to prove Theorem~\ref{app-mainthm}, we will first construct an equivalent game model whose states are the state of the original model plus the one-stage costs incurred up to that time. For instance, for the infinite-population limit, the new state is
\begin{align}
s(t) \coloneqq  \biggl(x(t),\sum_{k=0}^{t-1} c(x(k),a(k),\mu_k)\biggr). \nonumber
\end{align}
In this new model, finite-horizon risk-sensitive cost function can be written in an additive-form like in risk-neutral case. Therefore, we can use the proof technique in our previous paper \cite{SaBaRa17} to show the existence of an approximate Nash equilibrium. It is important to note that for this new game model, we have been inspired by \cite{BaRi14}, in which the authors study the classical risk-sensitive control problem. This new game model is specified by
\begin{align}
\biggl( \sS, \sA, \{P_t\}_{t\geq0}, \{C_t\}_{t\geq0}, \kappa_0 \biggr), \nonumber
\end{align}
where
\begin{align}
\sS &= \sX \times [0,L] \nonumber
\end{align}
with $L \coloneqq \frac{K}{1-\beta}$ and $\sA$ are the Polish state and action spaces, respectively. The stochastic kernel $P_t : \sS \times \sA \times \P(\sS) \to \P(\sS)$ is defined as:
\begin{align}
P_t\bigl(B \times D \big| s(t),a(t),\Delta_t\bigr) \coloneqq p(B|x(t),a(t),\Delta_{t,1}) \otimes \delta_{c(t) + \beta^t c(x(t),a(t),\Delta_{t,1})}(D), \nonumber
\end{align}
where $B \in \B(\sX)$, $D \in \B([0,L])$,
\begin{align}
s(t) = (x(t),c(t)), \nonumber
\end{align}
and $\Delta_{t,1}$ is the marginal of $\Delta_t$ on $\sX$. Here, $P_t$ is indeed the controlled transition probability of next state $x(t+1)$ and current  total cost $\sum_{k=0}^{t} \beta^k c(x(k),a(k),\Delta_{k,1})$ given the current state-action pair $(x(t),a(t))$ and past total cost $\sum_{k=0}^{t-1} \beta^k c(x(k),a(k),\Delta_{k,1})$ in the original mean-field game. For each $t\geq0$, the one-stage cost function $C_t: \sS \times \sA \times \P(\sS) \rightarrow [0,\infty)$ is defined as:
\begin{align}
C_t(s(t),a(t),\Delta_t) \coloneqq
\begin{cases}
0, & \text{ } \text{ if $t \leq n$} \\
e^{\lambda c(t)}, & \text{ } \text{ if $t = n + 1$}.
\end{cases}\nonumber
\end{align}
Finally, the initial measure $\kappa_0$ is given by $\kappa_0(ds(0)) \coloneqq \mu_0(dx(0)) \otimes \delta_0(dc(0))$. Note that, in this equivalent game model, the finite-horizon is $n+1$ instead of $n$.

Suppose that Assumptions~\ref{as1} and \ref{as2} hold.
Then, for each $t\geq0$, the following are satisfied:
\begin{itemize}
\item [(I)] The one-stage cost function $C_t$ is bounded and continuous.
\item [(II)] The stochastic kernel $P_t$ is weakly continuous.
\end{itemize}
It is straightforward to prove that (I) and (II) hold since $c$ is continuous and $p$ is weakly continuous.

Recall the set of Markov policies $\sM$ in the original mean-field game. Note that $\sM$ is a subset of the set of Markov policies in the new model. For any measure flow ${\boldsymbol \Delta} = (\Delta_t)_{t\geq0}$, where $\Delta_t \in \P(\sS)$, we denote by $\hat{J}_{{\boldsymbol \Delta}}^{n}(\pi)$ the finite-horizon risk-sensitive cost of the policy $\pi \in \sM$ in this new model.

We also define the corresponding $N$ agent game as follows. We have the Polish state space $\sS$ and action space $\sA$. For every $t \in \{0,1,2,\ldots\}$ and every $i \in \{1,2,\ldots,N\}$, let $s^N_i(t) = (x^N_i(t),c^N_i(t)) \in \sS$ and $a^N_i(t) \in \sA$ denote the state and the action of Agent~$i$ at time $t$, and let
\begin{align}
\Delta_t^{(N)}(\,\cdot\,) \coloneqq \frac{1}{N} \sum_{i=1}^N \delta_{s_i^N(t)}(\,\cdot\,) \in \P(\sS) \nonumber
\end{align}
denote the empirical distribution of the state configuration at time $t$. The initial states $s^N_i(0)$ are independent and identically distributed according to $\kappa_0$, and, for each $t \ge 0$, the next-state configuration $(s^N_1(t+1),\ldots,s^N_N(t+1))$ is generated at random according to the probability laws
\begin{align}
&\prod^N_{i=1} P_{t}\big(ds^N_i(t+1)\big|s^N_i(t),a^N_i(t),\Delta^{(N)}_t\big). \nonumber 
\end{align}

Recall that $\sM_i$ denotes the set of Markov policies for Agent $i$ in the original game. Note that $\sM_i$ is a subset of the Markov policies of the new model. Recall that $\sM_i^c$ denotes the set of all weakly continuous Markov policies for Agent~$i$ in the original game. For Agent~$i$, the finite-horizon risk-sensitive cost under the initial distribution $\kappa_0$ and $N$-tuple of policies ${\boldsymbol \pi}^{(N)} \in {\bf \sM}^{(N)}$ is denoted as $\hat{J}_i^{(N),n}({\boldsymbol \pi}^{(N)})$. The following proposition makes the connection between this new model and the original model.

\begin{proposition}\label{app-prop1}
For any $N\geq1$, ${\boldsymbol \pi}^{(N)} \in {\bf \sM}^{(N)}$, and $i=1,\ldots,N$, we have $\hat{J}_i^{(N),n}({\boldsymbol \pi}^{(N)}) = J_i^{(N),n}({\boldsymbol \pi}^{(N)})$. Similarly, for any $\pi \in \sM$ and measure flow ${\boldsymbol \Delta}$, we have $\hat{J}^n_{\boldsymbol \Delta}(\pi) = J^n_{\bmu}(\pi)$ where $\bmu = (\Delta_{t,1})_{t\geq0}$.
\end{proposition}

\proof{}
The proof of the proposition is given in Appendix~\ref{app1}.\Halmos
\endproof

By Proposition~\ref{app-prop1}, in the remainder of this section we consider the new game model in place of the original one. Then, we use a similar technique as in our previous paper \cite{SaBaRa17} to prove Theorem~\ref{app-mainthm}.

Define the measure flow ${\boldsymbol \Delta} = (\Delta_t)_{t\geq0}$ as follows: $\Delta_t = {\cal L}(x(t),\sum_{k=0}^{t-1} c(x(k),a(k),\mu_k))$, where ${\cal L}(x(t),\sum_{k=0}^{t-1} c(x(k),a(k),\mu_k))$ denotes the probability law of
$(x(t),\sum_{k=0}^{t-1} c(x(k),a(k),\mu_k))$ in the original mean-field game under the policy $\pi$ and measure flow $\bmu$ in the mean-field equilibrium. Note that, for each $t$, $\Delta_{t,1} = \mu_t$.
Define the stochastic kernel $P_t^{\pi}(\,\cdot\,|s,\Delta)$ on $\sS$ given $\sS \times \P(\sS)$ as
\begin{align}
P_t^{\pi}(\,\cdot\,|s,\Delta) \coloneqq \int_{\sA} P_t(\,\cdot\,|s,a,\Delta) \pi_t(da|s). \nonumber
\end{align}
Since $\pi_t$ is assumed to be weakly continuous, $P_t^{\pi}(\,\cdot\,|s,\Delta)$ is also weakly continuous in $(s,\Delta)$. In the sequel, to ease the notation, we will also write $P_t^{\pi}(\,\cdot\,|s,\Delta)$ as $P_{t,\Delta}^{\pi}(\,\cdot\,|s)$. Note that $\Delta_0 = \kappa_0$.

\begin{lemma}\label{app-lemma1}
Measure flow ${\boldsymbol \Delta}$ satisfies
\begin{align}
\Delta_{t+1}(\,\cdot\,) &= \int_{\sS} P_{t}^{\pi}(\,\cdot\,|s,\Delta_t) \Delta_t(ds) \nonumber \\
&= \Delta_t P_{t,\Delta_t}^{\pi}(\,\cdot\,). \nonumber
\end{align}
\end{lemma}

\proof{}

For any $t$, we have
\begin{align}
\Delta_{t+1}(\cdot) &= \int_{\sX \times \sA} p(\cdot|x,a,\mu_t) \otimes \delta_{c+\beta^t c(x,a,\mu_t)}(\cdot)\pi_t(da|x) \Delta_t(dx,dc) \nonumber \\
&= \int_{\sS \times \sA} P_t(\cdot|s,a,\Delta_t) \pi_t(da|x) \Delta_t(ds) \text{ } \text{ (as $\Delta_{t,1} = \mu_t$)}\nonumber \\
&= \Delta_t P_{t,\Delta_t}^{\pi}(\,\cdot\,), \nonumber
\end{align}
where $s = (x,c)$. This completes the proof.\Halmos
\endproof

For each $N\geq1$, let $\bigl\{s_i^{N}(t)\bigr\}_{1\leq i\leq N}$ denote the state configuration at time $t$ in the $N$-person game under the policy ${\boldsymbol \pi}^{(N)} = \{\pi,\pi,\ldots,\pi\}$. Define the empirical distribution
\begin{align}
\Delta_t^{(N)}(\,\cdot\,) \coloneqq \frac{1}{N} \sum_{i=1}^N \delta_{s_i^{N}(t)}(\,\cdot\,). \nonumber
\end{align}

\begin{proposition}\label{prop5}
For all $t\geq0$, we have
\begin{align}
{\cal L}(\Delta_t^{(N)}) \rightarrow \delta_{\Delta_t} \nonumber
\end{align}
weakly in $\P(\P(\sS))$, as $N\rightarrow\infty$.
\end{proposition}

\proof{}
The proof of the proposition is given in Appendix~\ref{app2}.\Halmos
\endproof

Proposition~\ref{prop5} states that, in the infinite-population limit, the empirical distribution of the states under the mean-field policy converges to the deterministic measure flow ${\boldsymbol \Delta}$. Since $P_t(\,\cdot\,|s,a,\Delta)$ is continuous in $\Delta$ for each $t$, the evolution of the state of a generic agent in the finite-agent game with sufficiently many agents and the evolution of the state in the mean-field game under policies $\bpi^{(N)} = (\pi,\ldots,\pi)$ and $\pi$, respectively, should therefore be close. Hence, the distributions of the states in each problem should also be close, from which we obtain the first part of Theorem~\ref{app-mainthm} in view of Proposition~\ref{app-prop1}.

\begin{proposition}\label{prop6}
We have
\begin{align}
\lim_{N\rightarrow\infty} \hat{J}_1^{(N),n}({\boldsymbol \pi}^{(N)}) = \hat{J}^{n}_{{\boldsymbol \Delta}}(\pi),  \nonumber
\end{align}
where $\bpi^{(N)} = (\pi,\ldots,\pi)$.
\end{proposition}

\proof{}
Recall that $C_t = 0$ for $t \leq n$. Hence, we let $t = n+1$. Since for any fixed permutation $\sigma$ of $\{1,\ldots,N\}$, we have
\begin{align}
{\cal L}\bigl(s_1^N(t),\ldots,s_N^N(t),\Delta_t^{(N)}\bigr) = {\cal L}\bigl(s_{\sigma(1)}^N(t),\ldots,s_{\sigma(N)}^N(t),\Delta_t^{(N)}\bigr), \nonumber
\end{align}
the cost function at time $t = n+1$ can be written as
\begin{align}
E\bigl[ C_{t}(s_1^N(t)) \bigr] &= \frac{1}{N} \sum_{i=1}^N E\bigl[ C_{t}(s_i^N(t)) \bigr] \nonumber \\
&= E\bigl[ \Delta_t^{(N)}\bigl(C_{t}(s)\bigr) \bigr]. \nonumber
\end{align}
Let $F: \P(\sS) \rightarrow \R$ be defined as
\begin{align}
F(\Delta) \coloneqq \int_{\sS} C_{t}(s) \Delta(ds). \nonumber
\end{align}
Note that $F \in C_b(\P(\sS))$ as $C_t \in C_b(\sS)$. Hence, by Proposition~\ref{prop5}, we obtain
\begin{align}
\lim_{N\rightarrow\infty} E\bigl[ C_{t}(s_1^N(t)) \bigr] &= \lim_{N\rightarrow\infty} E\bigl[ \Delta_t^{(N)}\bigl(C_{t}(s)\bigr) \bigr] \nonumber \\
&= \lim_{N\rightarrow\infty} E[F(\Delta_t^{(N)})] \nonumber \\
&= F(\Delta_t) \nonumber \\
&= \Delta_t(C_{t}). \label{eq9}
\end{align}
Note that by Lemma~\ref{app-lemma1}, the discounted cost in the mean-field game can be written as
\begin{align}
\hat{J}^n_{{\boldsymbol \Delta}}(\pi) = \Delta_t(C_{t}). \nonumber
\end{align}
Therefore, by (\ref{eq9}), we obtain
\begin{align}
\lim_{N\rightarrow\infty} \hat{J}_1^{(N),n}({\boldsymbol \pi}^{(N)}) = \hat{J}^n_{{\boldsymbol \Delta}}(\pi), \nonumber
\end{align}
which completes the proof.\Halmos
\endproof

Recall the policies $\{\tpi^{(N)}\}_{N\geq1} \subset \sM_1^c$ in Theorem~\ref{app-mainthm}. For each $N\geq1$, let $\bigl\{\ts_i^N(t)\bigr\}_{1\leq i \leq N}$ be the collection of states in the $N$-person game under the policy $\tilde{{\boldsymbol \pi}}^{(N)} \coloneqq \{\tpi^{(N)},\pi,\ldots,\pi\}$. Define
\begin{align}
\tilde{\Delta}_t^{(N)}(\,\cdot\,) \coloneqq \frac{1}{N} \sum_{i=1}^N \delta_{\ts_i^{(N)}(t)}(\,\cdot\,). \nonumber
\end{align}
The following result states that, in the infinite-population limit, the law of the empirical distribution of the states at each time $t$ is insensitive to local deviations from the mean-field equilibrium policy.

\begin{proposition}\label{prop8}
For all $t\geq0$, we have
\begin{align}
{\cal L}(\tilde{\Delta}_t^{(N)}) \rightarrow \delta_{\Delta_t} \nonumber
\end{align}
weakly $\P(\P(\sS))$, as $N \rightarrow \infty$.
\end{proposition}

\proof{}
The proof can be done by slightly modifying the proof of Proposition~\ref{prop5}, and therefore will not be included here. See the proof of \cite[Proposition 4.6]{SaBaRa17}.\Halmos
\endproof

For each $N\geq1$, let $\{\hs^N(t)\}_{t\geq0}$ denote the state trajectory of the mean-field game under policy $\tpi^{(N)}$; that is, $\hs^N(t)$ evolves as follows:
\begin{align}
\hs^N(0) \sim \kappa_0 \text{ and } \hs^N(t+1) \sim P^{\tpi^{(N)}}_{t,\Delta_t}(\,\cdot\,|\hs^N(t)). \nonumber
\end{align}
Recall that the cost function of this mean-field game is given by
\begin{align}
\hat{J}^n_{{\boldsymbol \Delta}}(\tpi^{(N)}) = E\bigl[ C_{n+1}(\hs^N(n+1))\bigr]. \label{nneq1}
\end{align}

The following result is very important for proving the second part of Theorem~\ref{app-mainthm}.

\begin{proposition}\label{prop9}
For any $t\geq0$, we have
\begin{align}
\lim_{N\rightarrow\infty} \bigl| {\cal L}(\ts_1^N(t))(g_N) - {\cal L}(\hs^N(t))(g_N) \bigr| = 0 \nonumber
\end{align}
for any sequence $\{g_N\} \subset C_b(\sS)$ such that $\sup_{N\geq1}\|g_N\|<\infty$ and $\omega_g(r) \rightarrow 0$ as $r \rightarrow 0$, where
\begin{align}
\omega_g(r) \coloneqq \sup_{x \in \sX} \sup_{N\geq1} \sup_{\substack{c,c' \\ |c - c'| \leq r}} |g_N(x,c) - g_N(x,c')|. \nonumber
\end{align}
\end{proposition}

\proof{}
The proof of the proposition is given in Appendix~\ref{app3}.\Halmos
\endproof

Using Proposition~\ref{prop9}, we now prove the following result, from which the second part of Theorem~\ref{app-mainthm} follows in view of Proposition~\ref{app-prop1}.

\begin{proposition}\label{theorem3}
Let $\{\tpi^{(N)}\}_{N\geq1} \subset \sM_1^c$ be an arbitrary sequence of policies for Agent~$1$. Then, we have
\begin{align}
\lim_{N \rightarrow \infty} \bigl| \hat{J}_1^{(N),n}(\tpi^{(N)},\pi,\ldots,\pi) - \hat{J}^n_{{\boldsymbol \Delta}}(\tpi^{(N)}) \bigr| = 0, \nonumber
\end{align}
where $\hat{J}^n_{{\boldsymbol \Delta}}(\tpi^{(N)})$ is given in (\ref{nneq1}).
\end{proposition}

\proof{}

Since $C_t = 0$ for $t \leq n$, we set $t=n+1$. We have
\begin{align}
\bigl| \hat{J}_1^{(N),n}(\tpi^{(N)},\pi,\ldots,\pi) - \hat{J}^n_{{\boldsymbol \Delta}}(\tpi^{(N)}) \bigr|
= \bigl| E\bigl[ C_t(\ts_1^N(t)) \bigr] - E\bigl[ C_t(\hs_1^N(t)) \bigr] \bigr|. \nonumber
\end{align}
Note that $C_t(s) = C_t((x,c)) \coloneqq e^{\lambda c}$, where $c \in [0,L]$, is Lipschitz function. Therefore, the term in the above equation converges to zero by Proposition~\ref{prop9}.\Halmos
\endproof

\proof{Proof of Theorem~\ref{app-mainthm}.}
The proof follows from Propositions~\ref{app-prop1}, \ref{prop6}, and \ref{theorem3}.\Halmos
\endproof

As a corollary of Lemma~\ref{lemma3} and Theorem~\ref{app-mainthm}, we obtain the following result.

\begin{corollary}\label{cor1}
We have
\begin{align}
\lim_{N \rightarrow \infty} J_1^{(N)}(\tpi^{(N)},\pi,\ldots,\pi)
&\geq \inf_{\pi' \in \sM} J_{{\boldsymbol \mu}}(\pi') \nonumber \\
&= J_{{\boldsymbol \mu}}(\pi) \nonumber \\
&= \lim_{N \rightarrow \infty} J_1^{(N)}(\pi,\pi,\ldots,\pi). \nonumber
\end{align}
\end{corollary}

\proof{}
Since, by Lemma~\ref{lemma3}, any infinite-horizon risk sensitive cost can be approximated by finite-horizon risk sensitive cost functions with error bound $L_0 \beta^{n+1}$, the result follows from Theorem~\ref{app-mainthm}.\Halmos
\endproof

Now, we are ready to prove the main result of this section.

\proof{Proof of Theorem~\ref{appr-thm}.}
One can prove that, for any policy ${\boldsymbol \pi}^{(N)} \in {\bf \sM}^{(N)}$, we have
\begin{align}
\inf_{\pi^i \in \sM_i} J_i^{(N)}({\boldsymbol \pi}^{(N)}_{-i},\pi^i) = \inf_{\pi^i \in \sM_i^c} J_i^{(N)}({\boldsymbol \pi}^{(N)}_{-i},\pi^i) \nonumber
\end{align}
for each $i=1,\ldots,N$ (see the proof of \cite[Theorem 2.3]{SaBaRa17}). Hence, it is sufficient to consider weakly continuous policies in ${\bf \sM}^{(N)}$ to establish the existence of $\varepsilon$-Nash equilibrium in the new model.

We should prove that, for sufficiently large $N$, we have
\begin{align}
J_i^{(N)}({\boldsymbol \pi}^{(N)}) &\leq \inf_{\pi^i \in \sM_i^c} J_i^{(N)}({\boldsymbol \pi}^{(N)}_{-i},\pi^i) + \varepsilon \label{eq13}
\end{align}
for each $i=1,\ldots,N$. As indicated earlier, since the transition probabilities and the one-stage cost functions are the same for all agents in the new game, it is sufficient to prove (\ref{eq13}) for Agent~$1$ only. Given $\epsilon > 0$, for each $N\geq1$, let $\tpi^{(N)} \in \sM_1^c$ be such that
\begin{align}
J_1^{(N)} (\tpi^{(N)},\pi,\ldots,\pi) < \inf_{\pi' \in \sM_1^c} J_1^{(N)} (\pi',\pi,\ldots,\pi) + \frac{\varepsilon}{3}. \nonumber
\end{align}
Then, by Corollary~\ref{cor1}, we have
\begin{align}
\lim_{N\rightarrow\infty} J_1^{(N)} (\tpi^{(N)},\pi,\ldots,\pi) &= \lim_{N\rightarrow\infty} J_{{\boldsymbol \mu}}(\tpi^{(N)}) \nonumber \\
&\geq \inf_{\pi'} J_{{\boldsymbol \mu}}(\pi') \nonumber \\
&= J_{{\boldsymbol \mu}}(\pi) \nonumber \\
&= \lim_{N\rightarrow\infty} J_1^{(N)} (\pi,\pi,\ldots,\pi). \nonumber
\end{align}
Therefore, there exists $N(\varepsilon)$ such that for $N\geq N(\varepsilon)$, we have
\begin{align}
\inf_{\pi' \in \sM_1^c} J_1^{(N)} (\pi',\pi,\ldots,\pi) + \varepsilon &> J_1^{(N)} (\tpi^{(N)},\pi,\ldots,\pi) + \frac{2\varepsilon}{3} \nonumber \\
&\geq J_{{\boldsymbol \mu}}(\pi) + \frac{\varepsilon}{3} \nonumber \\
&\geq J_1^{(N)} (\pi,\pi,\ldots,\pi). \nonumber
\end{align}
\Halmos
\endproof

\section{Conclusion}\label{conc}

This paper has studied discrete-time risk-sensitive stochastic games of the mean-field type. Letting the number of agents go to infinity, we have first established the existence of a mean-field equilibrium in the limiting mean-field game problem. We have then showed that the policy in the mean-field equilibrium constitutes an approximate Nash equilibrium for similarly structured games with a sufficiently large number of agents.

One interesting future direction would be to study risk-sensitive mean-field games with imperfect information, as the counterpart of the risk-neutral game studied in \cite{SaBaRa18}. In this case, a possible approach is to use the theory developed for partially observed risk-sensitive Markov decision processes. Finally, establishing similar results for different or more general risk-sensitivity criteria is also worth studying. In particular, mean-field games with conditional risk measures is an interesting research direction.

\section*{Appendix}

\appsec

\subsection{Proof of Proposition~\ref{app-prop1}}\label{app1}

Fix any $N\geq1$ and ${\boldsymbol \pi}^{(N)} \in {\bf \sM}^{(N)}$. For each $t\geq0$, let $\hat{\rP}_t$
denote the probability law of the states $\bs^N(t) = (s^N_1(t),\ldots,s^N_N(t))$ and the actions $\ba^N(t) = (a^N_1(t),\ldots,a^N_N(t))$ under ${\boldsymbol \pi}^{(N)}$ in the new finite-agent game model. Similarly, let $\rP_t$
denote the probability law of the states $\bx^N(t) = (x^N_1(t),\ldots,x^N_N(t))$, the total costs $\bc^N(t) = \biggl(\sum_{k=0}^{t-1} c(x_1^N(k),a_1^N(k),e_k^{(N)}),\ldots,\sum_{k=0}^{t-1} c(x_N^N(k),a_N^N(k),e_k^{(N)})\biggr)$, and the actions $\ba^N(t) = (a^N_1(t),\ldots,a^N_N(t))$ under ${\boldsymbol \pi}^{(N)}$ in the original finite agent game model. Here, we assume that $\sum_{k=0}^{-1} c(x_i^N(k),a_i^N(k),e_k^{(N)}) = 0$ for any $i=1,\ldots,N$. We prove that, for each $t\geq0$,
\begin{align}
\hat{\rP}_t = \rP_t, \nonumber
\end{align}
which implies that $\hat{J}_i^{(N),n}({\boldsymbol \pi}^{(N)}) = J_i^{(N),n}({\boldsymbol \pi}^{(N)})$ for all $i=1,\ldots,N$.

The claim trivially holds for $t=0$. Suppose that the claim holds for $t$ and consider $t+1$. For $t+1$, let
\begin{align}
{\boldsymbol A}_{t+1} &= A_{t+1}^1 \times \ldots \times A_{t+1}^N, \text{ } \text{$A_{t+1}^i \in \B(\sX)$ for $i=1,\ldots,N$}, \nonumber \\
{\boldsymbol B}_{t+1} &= B_{t+1}^1 \times \ldots \times B_{t+1}^N, \text{ } \text{$B_{t+1}^i \in \B([0,L])$ for $i=1,\ldots,N$}, \nonumber \\
{\boldsymbol D}_{t+1} &= D_{t+1}^1 \times \ldots \times D_{t+1}^N, \text{ } \text{$D_{t+1}^i \in \B(\sA)$ for $i=1,\ldots,N$}. \nonumber
\end{align}
Define $G_{t+1} \coloneqq {\boldsymbol A}_{t+1} \times {\boldsymbol B}_{t+1} \times {\boldsymbol D}_{t+1}$. Then, we have
\begin{align}
&\hat{\rP}_t(G_{t+1}) \nonumber \\
&= \int_{\sS \times \sA} \int_{G_{t+1}} \biggl( \prod_{i=1}^N \pi_{t+1}^i(da^N_i(t+1)|s^N_i(t+1))  \nonumber \\
&\phantom{xxxxxxxxxxxx} P_t(ds^N_i(t+1)|s^N_i(t),a_i^N(t),\Delta_t^{(N)}) \biggr) \hat{\rP}_t(d\bs^N(t),d\ba^N(t))  \nonumber \\
&= \int_{\sS \times \sA} \int_{G_{t+1}} \biggl( \prod_{i=1}^N \pi_{t+1}^i(da^N_i(t+1)|x^N_i(t+1))  \nonumber \\ &\phantom{xx}P_t(dx^N_i(t+1),dc_i^N(t+1)|x^N_i(t),c_i^N(t),a_i^N(t),\Delta_t^{(N)}) \biggr) \hat{\rP}_t(d\bs^N(t),d\ba^N(t))  \nonumber \\
&= \int_{\sS \times \sA} \int_{G_{t+1}}  \biggl( \prod_{i=1}^N \pi_{t+1}^i(da^N_i(t+1)|x^N_i(t+1))  \nonumber \\
&\phantom{xxxx} p(dx^N_i(t+1)|x^N_i(t),a^N_i(t),e_t^{(N)}) \delta_{c_i^N(t) + \beta^t c(x_i^N(t),a_i^N(t),e_t^{(N)})}(dc_i^N(t+1)) \biggr) \nonumber \\
&\phantom{xxxxxxxxxxxx}\rP_t(d\bx^N(t),d\bc^N(t),d\ba^N(t)) \text{ } \text{ (by induction hypothesis)}  \nonumber \\
&= \rP_{t+1}(G_{t+1}). \nonumber
\end{align}
Hence, $\hat{\rP}_{t+1} = \rP_{t+1}$. This implies that $\hat{J}_i^{(N),n}({\boldsymbol \pi}^{(N)}) = J_i^{(N),n}({\boldsymbol \pi}^{(N)})$ for all $i=1,\ldots,N$.

The second part of the proposition can be proved similarly, so we omit the details.

\subsection{Proof of Proposition~\ref{prop5}}\label{app2}

Although this result can be deduced from \cite[Proposition 4.4]{SaBaRa17}, we give the proof in detail for the sake of completeness.

It is known that weak topology on $\P(\sS)$ can be metrized using the following metric:
\begin{align}
\rho(\mu,\nu) \coloneqq \sum_{m=1}^{\infty} 2^{-(m+1)} | \mu(f_m) - \nu(f_m) |, \nonumber
\end{align}
where $\{f_m\}_{m\geq1}$ is an appropriate sequence of real-valued continuous and bounded functions on $\sS$ such that $\|f_m\| \leq 1$ for all $m\geq1$ (see \cite[Theorem 6.6, p. 47]{Par67}). Define the Wasserstein distance of order 1 on the set of probability measures $\P(\P(\sS))$ as follows (see \cite[Definition 6.1]{Vil09}):
\begin{align}
W_1(\Phi,\Psi) \coloneqq \inf \bigl\{ E[\rho(X,Y)]: {\cal L}(X) = \Phi \text{ and } {\cal L}(Y) = \Psi \bigr\}. \nonumber
\end{align}
Note that since $\delta_{\Delta_t}$ is a Dirac measure, we have
\begin{align}
W_1({\cal L}(\Delta_t^{(N)}),\delta_{\Delta_t}) &= \bigl\{ E[\rho(X,Y)]: {\cal L}(X) = {\cal L}(\Delta_t^{(N)}) \text{ and } {\cal L}(Y) = \delta_{\Delta_t} \bigr\} \nonumber \\
&= E\biggl[ \sum_{m=1}^{\infty} 2^{-(m+1)} | \Delta_t^{(N)}(f_m) - \Delta_t(f_m) | \biggr]. \nonumber
\end{align}
Since convergence in $W_1$ distance implies weak convergence (see \cite[Theorem 6.9]{Vil09}), it suffices to prove that
\begin{align}
\lim_{N\rightarrow\infty} E\bigl[|\Delta_t^{(N)}(f) - \Delta_t(f)|\bigr] = 0 \nonumber
\end{align}
for any $f \in C_b(\sS)$ and for all $t$. We prove this by induction on $t$.

As $\{s_i^N(0)\}_{1\leq i\leq N} \sim \kappa_0^{\otimes N}$, the claim is true for $t=0$. We suppose that the claim holds for $t$ and consider $t+1$. Fix any $g \in C_b(\sS)$. Then, we have
\begin{align}
&|\Delta_{t+1}^{(N)}(g) - \Delta_{t+1}(g)| \leq |\Delta_{t+1}^{(N)}(g) - \Delta_{t}^{(N)} P^{\pi}_{t,\Delta_t^{(N)}}(g)| \nonumber \\
&\phantom{xxxxxxxxxxxxxxxxxxxxxxxxx} + |\Delta_t^{(N)} P^{\pi}_{t,\Delta_t^{(N)}}(g) - \Delta_t P^{\pi}_{t,\Delta_t}(g) |. \label{eq8}
\end{align}
We first prove that the expectation of the second term on the right-hand side (RHS) of (\ref{eq8}) converges to $0$ as $N\rightarrow\infty$. To that end, define $F: \P(\sS) \rightarrow \R$ as
\begin{align}
F(\Delta) = \Delta P^{\pi}_{t,\Delta}(g) \coloneqq \int_{\sS} \int_{\sS} g(s') P^{\pi}_t(ds'|s,\Delta) \Delta(ds). \nonumber
\end{align}
One can prove that $F \in C_b(\P(\sS))$. Indeed, suppose that $\Delta_n$ converges to $\Delta$. Let us define
\begin{align}
l_n(s) &= \int_{\sS} g(s') P^{\pi}_t(ds'|s,\Delta_n) \nonumber \\
\intertext{and}
l(s) &= \int_{\sS} g(s') P^{\pi}_t(ds'|s,\Delta). \nonumber
\end{align}
Since $P^{\pi}_t$ is weakly continuous, one can prove that $l_n$ converges to $l$ continuously; that is, if $s_n$ converges to $s$, then $l_n(s_n) \rightarrow l(s)$. By \cite[Theorem 3.5]{Lan81}, we have $F(\Delta_n) \rightarrow F(\Delta)$, and so, $F \in C_b(\P(\sS))$. This implies that the expectation of the second term on the RHS of (\ref{eq8}) converges to zero as ${\cal L}(\Delta_t^{(N)}) \rightarrow \Delta_t$ weakly, by the induction hypothesis.

Now, let us write the expectation of the first term on the RHS of (\ref{eq8}) as
\begin{align}
E\biggl[ E\biggl[ |\Delta_{t+1}^{(N)}(g) - \Delta_t^{(N)} P^{\pi}_{t,\Delta_t^{(N)}}(g)| \biggr| s_1^N(t),\ldots,s_N^N(t) \biggr] \biggr]. \nonumber
\end{align}
Then, by \cite[Lemma A.2]{BuMa14}, we have
\begin{align}
E\biggl[ |\Delta_{t+1}^{(N)}(g) - \Delta_t^{(N)} P^{\pi}_{t,\Delta_t^{(N)}}(g)| \biggr| s_1^N(t),\ldots,s_N^N(t) \biggr] \leq 2 \frac{\|g\|}{\sqrt{N}}. \nonumber
\end{align}
Therefore, the expectation of the first term on the RHS of (\ref{eq8}) also converges to zero as $N\rightarrow\infty$. Since $g$ is arbitrary, this completes the proof.

\subsection{Proof of Proposition~\ref{prop9}}\label{app3}

Let $\{T_N\} \subset C_b(\sS\times\P(\sS))$ be an arbitrary sequence of functions such that the family $\bigl\{T_N(s,\,\cdot\,): s \in \sS, N\geq1)\bigr\}$ is equi-continuous, $\sup_{N\geq1}\|T_N\|<\infty$, and $\omega_{T,\Delta}(r) \rightarrow 0$ as $r \rightarrow 0$ for any $\Delta \in \P(\sS)$, where
\begin{align}
\omega_{T,\Delta}(r) \coloneqq \sup_{x \in \sX} \sup_{N\geq1} \sup_{\substack{c,c' \\ |c - c'| \leq r}} |T_N(x,c,\Delta) - T_N(x,c',\Delta)|. \nonumber
\end{align}

Fix any $t\geq0$ and suppose that the statement of the proposition holds for this $t$; that is, we have
\begin{align}
\lim_{N\rightarrow\infty}\bigl| {\cal L}(\ts_1^N(t))(g_N) - {\cal L}(\hs^N(t))(g_N) \bigr| = 0 \label{neweq1}
\end{align}
for any sequence $\{g_N\}_{N\geq1} \subset C_b(\sS)$ satisfying the hypothesis of the proposition. Given (\ref{neweq1}), we prove that
\begin{align}
\lim_{N\rightarrow\infty} \bigl| {\cal L}(\ts_1^N(t),\tilde{\Delta}_t^{(N)})(T_N) - {\cal L}(\hs^N(t),\delta_{\Delta_t})(T_N) \bigr| = 0. \label{neweq2}
\end{align}
Indeed, we have
\begin{align}
&\bigl| {\cal L}(\ts_1^N(t),\tilde{\Delta}_t^{(N)})(T_N) - {\cal L}(\hs^N(t),\delta_{\Delta_t})(T_N) \bigr| \nonumber \\
&\phantom{xxxxxx}\leq \biggl| \int_{\sS \times \P(\sS)} T_N(s,\Delta) {\cal L}(\ts_1^N(t),\tilde{\Delta}_t^{(N)})(ds,d\Delta) \nonumber \\
&\phantom{xxxxxxxxxxxxxxxxxxx}- \int_{\sS \times \P(\sS)} T_N(s,\Delta) {\cal L}(\ts_1^N(t),\delta_{\Delta_t})(ds,d\Delta) \biggr| \nonumber \\
&\phantom{xxxxxx}+ \biggl| \int_{\sS \times \P(\sS)} T_N(s,\Delta) {\cal L}(\ts_1^N(t),\delta_{\Delta_t})(ds,d\Delta)  \nonumber \\
&\phantom{xxxxxxxxxxxxxxxxxxx}- \int_{\sS \times \P(\sS)} T_N(s,\Delta) {\cal L}(\hs^N(t),\delta_{\Delta_t})(ds,d\Delta) \biggr|. \label{eq12}
\end{align}
First, note that since the family $\{T_N(\,\cdot\,,\Delta_t)\}_{N\geq1} \subset C_b(\sS)$ satisfies the hypothesis of the proposition, we have
\begin{align}
\lim_{N\rightarrow\infty} \biggl| \int_{\sS} T_N(s,\Delta_t) {\cal L}(\ts_1^N(t))(ds) - \int_{\sS} T_N(s,\Delta_t) {\cal L}(\hs^N(t))(ds) \biggr|=0 \nonumber
\end{align}
by (\ref{neweq1}). Hence, the second term in (\ref{eq12}) converges to zero as $N\rightarrow\infty$.

Now, let us consider the first term in (\ref{eq12}). To that end, define ${\cal F} \coloneqq \bigl\{T_N(s,\,\cdot\,): s \in \sS, N\geq1)\bigr\}$. Note that ${\cal F}$ is a uniformly bounded and equi-continuous family of functions on $\P(\sS)$, and therefore
\begin{align}
\lim_{N\rightarrow\infty} E\biggl[ \sup_{F \in {\cal F}} \bigl| F(\tilde{\Delta}_t^{(N)}) - F(\Delta_t) \bigr| \biggr] = 0 \nonumber
\end{align}
as ${\cal L}(\tilde{\Delta}_t^{(N)}) \rightarrow {\cal L}(\Delta_t)$ weakly. Then, we have
\begin{align}
&\lim_{N\rightarrow\infty}\biggl| \int_{\sS \times \P(\sS)} T_N(s,\Delta) {\cal L}(\ts_1^N(t),\tilde{\Delta}_t^{(N)})(ds,d\Delta) \nonumber \\
&\phantom{xxxxxxxxxxxxxxxxxx}- \int_{\sS \times \P(\sS)} T_N(s,\Delta) {\cal L}(\ts_1^N(t),\delta_{\Delta_t})(ds,d\Delta) \biggr| \nonumber \\
&\leq \lim_{N\rightarrow\infty} \int_{\sS} \biggl| \int_{\P(\sS)} T_N(s,\Delta) {\cal L}(\tilde{\Delta}_t^{(N)}|\ts_1^N(t))(d\Delta|s) \nonumber \\
&\phantom{xxxxxxxxxxxxxxxxxx}- \int_{\P(\sS)} T_N(s,\Delta) {\cal L}(\delta_{\Delta_t})(d\Delta) \biggr| {\cal L}(\ts_1^N(t))(ds) \nonumber \\
&\leq \lim_{N\rightarrow\infty} E\biggl[ E\biggl[ \bigl|T_N(\ts_1^N(t),\tilde{\Delta}_t^{(N)}) - T_N(\ts_1^N(t),\Delta_t) \bigl| \biggl| \ts_1^N(t) \biggr] \biggr] \nonumber \\
&\leq \lim_{N\rightarrow\infty} E\biggl[ \sup_{F \in {\cal F}} \bigl|F(\tilde{\Delta}_t^{(N)}) - F(\Delta_t) \bigl| \biggr] \nonumber \\
&= 0. \nonumber
\end{align}
Hence, (\ref{neweq1}) implies (\ref{neweq2}) for any $t$.

Now, we prove that (\ref{neweq1}) is true for all $t$, which will complete the proof. Set $\sup_{N\geq1} \|g_N\| \eqqcolon L <\infty$ and define
\begin{align}
l_{N,t}(s,\Delta) \coloneqq \int_{\sA \times \sS} g_N(s') P_t(ds'|s,a,\Delta) \tpi_t^{(N)}(da|x), \nonumber
\end{align}
where $s = (x,c)$. For any $s \in \sS$ and $(\Delta,\Delta') \in \P(\sS)^2$, we have
\begin{align}
&|l_{N,t}(s,\Delta) - l_{N,t}(s,\Delta')| \nonumber \\
&\phantom{xxxxxxxx}= \biggl| \int_{\sA \times \sX} g_N(x',c+c(x,a,\Delta_1)) p(dx'|x,a,\Delta_1) \tpi_t^{(N)}(da|x) \nonumber \\
&\phantom{xxxxxxxxxxxxxx}- \int_{\sA \times \sX} g_N(x',c+c(x,a,\Delta'_1)) p(dx'|x,a,\Delta'_1) \tpi_t^{(N)}(da|x) \biggr| \nonumber \\
&\phantom{xxxxxxxx}\leq \biggl| \int_{\sA \times \sX} g_N(x',c+c(x,a,\Delta_1)) p(dx'|x,a,\Delta_1) \tpi_t^{(N)}(da|x) \nonumber \\
&\phantom{xxxxxxxxxxxxxx}- \int_{\sA \times \sX} g_N(x',c+c(x,a,\Delta'_1)) p(dx'|x,a,\Delta_1) \tpi_t^{(N)}(da|x) \biggr| \nonumber \\
&\phantom{xxxxxxxxxxx}+ \biggl| \int_{\sA \times \sX} g_N(x',c+c(x,a,\Delta'_1)) p(dx'|x,a,\Delta_1) \tpi_t^{(N)}(da|x) \nonumber \\
&\phantom{xxxxxxxxxxxxxx}- \int_{\sA \times \sX} g_N(x',c+c(x,a,\Delta'_1)) p(dx'|x,a,\Delta'_1) \tpi_t^{(N)}(da|x) \biggr| \nonumber \\
&\phantom{xxxxxxxx}\leq \int_{\sA} \omega_g(|c(x,a,\Delta_1) - c(x,a,\Delta'_1)|) \tpi_t^{(N)}(da|x) + L \omega_p(d_{BL}(\Delta_1,\Delta'_1)) \nonumber \\
&\phantom{xxxxxxxx}\leq \omega_g\bigl(\omega_c(d_{BL}(\Delta_1,\Delta'_1))\bigr) +  L \omega_p(d_{BL}(\Delta_1,\Delta'_1)). \nonumber
\end{align}
Hence the family $\{l_{N,t}(s,\,\cdot\,): s \in \sS, N\geq1\}$ is uniformly bounded and equi-continuous. Moreover, for any $\Delta \in \P(\sS)$, we have
\begin{align}
\omega_{l,\Delta}(r) &\coloneqq \sup_{x \in \sX} \sup_{N\geq1} \sup_{\substack{c,c' \\ |c - c'| \leq r}} |l_{N,t}(x,c,\Delta) - l_{N,t}(x,c',\Delta)| \nonumber \\
&= \sup_{x \in \sX} \sup_{N\geq1} \sup_{\substack{c,c' \\ |c - c'| \leq r}} \biggl| \int_{\sA \times \sX} g_N(x',c+c(x,a,\Delta_1)) p(dx'|x,a,\Delta_1) \tpi_t^{(N)}(da|x) \nonumber \\
&\phantom{xxxxxxxxxxxxxx}- \int_{\sA \times \sX} g_N(x',c'+c(x,a,\Delta_1)) p(dx'|x,a,\Delta_1) \tpi_t^{(N)}(da|x) \biggr| \nonumber \\
&\leq \sup_{x \in \sX} \sup_{N\geq1} \sup_{\substack{c,c' \\ |c - c'| \leq r}} \omega_g(|c-c'|) = \omega_g(r). \nonumber
\end{align}
Hence, $\omega_{l,\Delta}(r) \rightarrow 0$ as $r \rightarrow 0$.

We now prove (\ref{neweq1}) by induction on $t$. The claim trivially holds for $t=0$ as ${\cal L}(\ts_1^N(0)) = {\cal L}(\hs^N(0)) = \kappa_0$ for all $N\geq1$. Suppose the claim holds for $t$ and consider $t+1$. We can write
\begin{align}
&\bigl| {\cal L}(\ts_1^N(t+1))(g_N) - {\cal L}(\hs^N(t+1))(g_N) \bigr| \nonumber \\
&\phantom{xxxxx}=\biggl| \int_{\sS \times \P(\sS)} l_{N,t}(s,\Delta) {\cal L}(\ts_1^N(t),\tilde{\Delta}_t^{(N)})(ds,d\Delta) \nonumber \\
&\phantom{xxxxxxxxxxxxxxxxxxxxxxxx}- \int_{\sS \times \P(\sS)} l_{N,t}(s,\Delta) {\cal L}(\hs^N(t),\delta_{\Delta_t})(ds,d\Delta) \biggr|. \nonumber
\end{align}
Since the family $\{l_{N,t}\}_{N\geq1}$ is equi-continuous, uniformly bounded, $\omega_{l,\Delta}(r) \rightarrow 0$ as $r \rightarrow 0$ for any $\Delta \in \P(\sS)$ and (\ref{neweq1}) implies (\ref{neweq2}) at time $t$, the last term converges to zero as $N\rightarrow\infty$. This completes the proof.

\bibliographystyle{plain}

\end{document}